\documentclass[twoside,centertags]{amsart}

\newtheorem{innercustomthm}{Theorem}
\newenvironment{customthm}[1]
  {\renewcommand\theinnercustomthm{#1}\innercustomthm}
  {\endinnercustomthm}
\usepackage{mathtools}  
\usepackage{fancyhdr}            %
\fancypagestyle{plain}{
    \lhead{}
    \fancyhead[R]{\thepage}
    \fancyhead[L]{}
    
    \fancyfoot{}
}

\usepackage[cmtip,color,all]{xy}
\usepackage{amsmath}
\usepackage{amsfonts}
\usepackage{amssymb}
\usepackage{amsthm}
\usepackage{graphicx}
\usepackage{float}
\usepackage{needspace}
\usepackage{mathrsfs}
\usepackage[hyperfootnotes=false]{hyperref}
\usepackage{tikz-cd}
\usepackage{xcolor}
\usepackage{bbold}
\usepackage{pifont}
\usepackage{stmaryrd}
\SetSymbolFont{stmry}{bold}{U}{stmry}{m}{n}

\newtheorem{theorem}{Theorem}[subsection]
\newtheorem{corollary}[theorem]{Corollary}
\newtheorem{proposition}[theorem]{Proposition}

\theoremstyle{definition}
\newtheorem{definition}[theorem]{Definition}

\theoremstyle{definition}
\newtheorem{remark}[theorem]{Remark}

\usepackage{geometry}
\usepackage{setspace}
\usepackage{Commons,soul,tikz}
\usetikzlibrary{cd}
\usepackage{circuitikz}

\usepackage[backend=biber, style = alphabetic, doi = false, url=false, isbn = false]{biblatex}
\usepackage{xurl}
\hypersetup{breaklinks=true,
  pdftitle={Goodwillie polytopes: a family of self-dual prismatic polytopes},
  pdfauthor={Manuel Rivera and Daniel Tolosa}}
\renewbibmacro{in:}{}
\DeclareFieldFormat[article]{title}{\mkbibemph{#1}}
\DeclareFieldFormat[misc]{title}{\mkbibemph{#1}}
\DeclareFieldFormat{journaltitle}{#1}
\definecolor{winered}{rgb}{0.5,0,0}
 \hypersetup{citecolor = winered, linkcolor = winered, menucolor = winered, colorlinks = true}

\title{Goodwillie polytopes: a family of self-dual prismatic polytopes}
\author[M. Rivera]{Manuel Rivera}

\address{\newline M.R., Department of Mathematics, Purdue University, 150 N University St., West Lafayette, IN, USA 47907
\newline
Visiting Professor, Departamento de Matemáticas, Universidad Nacional de Colombia, Sede Bogotá, Carrera 45 No. 26-85, Bogotá D.C., Colombia}

\email{\href{mailto:manuelr@purdue.edu}{manuelr@purdue.edu}}

\author[D. Tolosa]{Daniel Tolosa}
\address{\newline D.T., School of Mathematical and Statistical Sciences, Arizona State University, 1151 S Forest Ave, Tempe, AZ, USA 85281}
\email{\href{mailto:dtolosav@asu.edu}{dtolosav@asu.edu}}

\begin{document}

\singlespacing

\begin{abstract} Families of self-dual polytopes with restricted facet types are hard to come by: for $n \geq 3$, the $n$-simplex is the only self-dual simplicial $n$-polytope, and there are no self-dual cubical $n$-polytopes. We prove that the Goodwillie polytopes form a family of self-dual prismatic polytopes, one in each dimension. For $n\geq3$, these do not arise as a join of two non-empty polytopes or as elongated pyramids (standard constructions for producing new self-dual polytopes from old ones). For each polytope in the family, we describe an explicit regular decomposition into prisms and give a closed formula for its $f$-polynomial. Finally, we show that Goodwillie polytopes arise as unions of prisms in a cell complex modeling the bar resolution of the path algebra of a simplicial complex: they are the images of generators of a smaller resolution under a quasi-isomorphism into the bar resolution.

\end{abstract}

\maketitle

\section{Introduction}

Beyond a few standard constructions, self-dual polytopes are hard to come by in higher dimensions, and self-dual families with restricted facet types are rarer still. In this article we study one such family, in which every facet of every polytope is a prism.  Our first main result is the following.

\begin{customthm}{1}\label{theorem1}
The Goodwillie polytopes $G_0,G_1,G_2,\ldots $ form a
family $\mathcal{G}$ of convex integral polytopes, where $\operatorname{dim}(G_n)=n$, with the following combinatorial properties.
\begin{enumerate}

  \item Each $G_n$ is combinatorially self-dual: its face poset is isomorphic to its dual face
        poset.
\item For $n\geq1$, $G_n$ is prismatic: every facet of $G_n$ is combinatorially equivalent to a prism $\Delta^k\times I^l$, the product of a $k$-simplex and an $l$-cube, where $k+l=n-1$.
\item For $n\ge 3$, $G_n$ is not combinatorially equivalent to a join of two non-empty polytopes or to an elongated pyramid.
\item The number $f_d(G_n)$ of $d$-dimensional faces of $G_n$ is a polynomially recursive function of $n$.
\end{enumerate}
\end{customthm}

 \textit{Goodwillie polytopes} were originally constructed by the authors while studying loop spaces and Koszul duality in \cite{rivera_cyclic_2024}. An explicit convex realization for the Goodwillie polytopes is given in Definition \ref{def:G_n}. 
 
The underlying combinatorial type of each $G_n$ may be obtained by appropriately gluing the elements in the set of prisms \[\{ \Delta^{|S|+1} \times I^{n-|S|-1} \}_{S \subseteq \{2,\ldots,n\}},\] as explained in Theorem \ref{thm:subdivision}. This decomposition is key in understanding the combinatorics of Goodwillie polytopes. 

A polytope is (combinatorially) \textit{self-dual} if its face poset is isomorphic to its dual face poset. Every $2$-dimensional polytope (a polygon) is self-dual. Simplices provide examples of self-dual polytopes in every dimension. Further examples can be constructed from existing ones through the following standard constructions. If $P$ is a self-dual $n$-dimensional polytope, then both the pyramid and elongated pyramid constructions applied to $P$ yield self-dual polytopes of dimension $n+1$. More generally, if $P$ and $Q$ are self-dual polytopes of dimensions $n$ and $m$, respectively, then their join $P \star Q$ is a self-dual polytope of dimension $n+m+1$. Thus, statement (3) can be interpreted as saying that Goodwillie polytopes are ``prime" self-dual polytopes in certain sense. 

There is an infinite number of combinatorially distinct self-dual convex polyhedra ($3$-dimensional polytopes) and all are classified through six possible constructions in \cite{AB92}. A good exercise is to examine these constructions carefully and classify all possible self-dual convex polyhedra with triangular and quadrilateral facets, i.e., self-dual prismatic convex polyhedra, see also \cite{Maf21}. Some examples of self-dual convex polyhedra are all pyramids and elongated pyramids over polygons, the diminished rhombic dodecahedron, the \textit{Josehedron} (a polyhedron recently described in \cite{MB26}), and $G_3$. Incidentally, the Josehedron has the same face-vector as $G_3$ (both have 12 vertices, 22 edges, and 8 quadrilateral faces and 4 triangular faces) but is not combinatorially equivalent to it. Both the Josehedron and the diminished rhombic dodecahedron happen to be prismatic and space-filling. An example of a $4$-dimensional self-dual polytope which is not a join of two polytopes or an elongated pyramid is the icositetrachoron. 

It is harder to encounter in ``nature" families of self-dual polytopes, or even individual higher dimensional examples, that are not given through the constructions described above. If one restricts the structure of the polytopes further, it is even harder to encounter self-duality. For instance, for $n \geq 3$, the only example, up to combinatorial equivalence, of a simplicial $n$-polytope (an $n$-polytope whose all facets are simplices) which is self-dual is the $n$-simplex. Therefore, there is no example of a family of polytopes satisfying (1), (3), and (4) with (2) replaced by the condition that each polytope is simplicial. Furthermore, there are no cubical $n$-polytopes ($n$-polytopes whose all facets are cubes) which are self-dual, for $n \geq 3$ \cite[Theorem 19.5.25]{Kal97}.  In the present article, we show that if, instead of the cubical or simplicial condition, we allow the polytopes to be \textit{prismatic} as described in (2), the Goodwillie polytopes provide an example of such a family. 

A function $f(n)$ is \textit{polynomially recursive} (or \textit{P-recursive}) if it satisfies a linear recurrence $p_0(n)f(n)+p_1(n)f(n+1)+\cdots+p_k(n)f(n+k)=0$ whose coefficients $p_i$ are polynomials in $n$, not all zero. In the case of Goodwillie polytopes, we prove something more explicit; for $n\geq3$ and $0\leq d\leq n-1$, we prove the following closed formula for the number $f_d(G_n)$ of $d$-dimensional faces of $G_n$: 
  \[f_d(G_n)
  =
  2^{n-1}\binom{n+1}{d+1}
  -2^{n-d-2}\binom{n-1}{d+1}
  -2^{d-1}\binom{n-1}{d-1}.
  \]

It follows that, for fixed $d$, $f_d(G_n)$  is a linear combination of $2^nn^i$ ($0\le i\le d+1$) and $n^j$ ($0\le j\le d-1$), so it satisfies a linear recurrence with \textit{constant} coefficients with characteristic polynomial $(x-2)^{d+2}(x-1)^d$. We state (4) in this form in order to pose the question about the possibility of having more families of polytopes that come from a unified combinatorial construction satisfying (1)-(4), see question (\ref{Q2}) below. The \textit{Ferrers polytopes}, as constructed in the PhD thesis \cite[Chapter 4]{St14}, is a family of polytopes satisfying (1), (3), and (4) with facets being products of simplices. 

What is most surprising about the family $\mathcal{G}$ is the context in which it was discovered: in \cite{rivera_cyclic_2024}, we discovered these polytopes while analyzing combinatorial models for loop spaces from the point of view of homological algebra. The name \textit{Goodwillie polytopes} came about since these were discovered after unraveling the combinatorics behind a map of chain complexes relating the Hochschild complex of the differential graded algebra of singular chains on the based loop space of a pointed topological space $X$ and the singular chains on the free loop space of $X$, originally constructed by Goodwillie in \cite{GOODWILLIE1985187}. 

In the second main result of this paper, stated as Theorem \ref{theorem2} below and discussed in a relatively self-contained manner in Section \ref{sec:resolutions}, we revisit this first encounter with Goodwillie polytopes from a purely combinatorial perspective. The structure of the Goodwillie polytopes reflects the combinatorics behind a resolution of certain differential graded algebra $\mathbb{A}(K)$ associated to any simplicial complex $K$, as described next. 

Let $K$ be a simplicial complex with a coherent ordering of the vertices of each simplex and $R$ any commutative ring. For any $n$-simplex $\sigma \in K$, denote by $s(\sigma)$ and $t(\sigma)$ the first and last vertices of $\sigma$, respectively. Denote by $C_0(K)$ the formal $R$-linear combinations on the vertices of $K$ considered as a cocommutative coalgebra with coproduct $\Delta \colon C_0(K) \to C_0(K) \otimes C_0(K)$ determined by declaring $\Delta(v)=v \otimes v$ for any vertex $v$ in $K$. 

We construct, directly from the combinatorics of $K$, a monoid object $\mathbb{A}(K)$ in the monoidal category of differential graded $C_0(K)$-bicomodules with cotensor product (denoted by $\square$). This monoid $\mathbb{A}(K)$, called the \textit{path algebra} of $K$, is a higher dimensional analog of the path algebra of a quiver. It encodes, in a single package, the quasi-isomorphism type of all chain complexes $S_*(P_{v,w}|K|;R)$ of singular $R$-chains in the space $P_{v,w}|K|$ of continuous paths between vertices $v$ and $w$ in the geometric realization of $K$. Roughly, $\mathbb{A}(K)$ is constructed as follows:
\begin{enumerate}
\item Consider the free graded $R$-module generated by \textit{necklaces} in $K$ (sequences $[\sigma_1|\cdots|\sigma_k]$ where each $\sigma_i$ is a positive-dimensional simplex in $K$ and $t(\sigma_i)= s(\sigma_{i+1})$) together with the symbols $\text{id}_v$, for all vertices $v$ in $K$. The degree of $[\sigma_1|\cdots|\sigma_k]$ is declared to be $\text{dim}(\sigma_1) + \cdots + \text{dim}(\sigma_k)-k$ and the degree of all $\text{id}_v$ is zero. 
\item On the resulting graded $R$-module, define a $C_0(K)$-bicomodule structure through the maps $[\sigma_1|\cdots|\sigma_k] \mapsto s(\sigma_1)$ and $[\sigma_1|\cdots|\sigma_k] \mapsto t(\sigma_k)$.
\item Consider the resulting $C_0(K)$-bicomodule as a free monoid with product given by concatenating necklaces.
\item In the resulting monoid, add formal inverses for all $[\gamma]$ where $\gamma$ is a $1$-simplex in $K$.
\item Define a differential by extending the sum of the simplicial boundary and the Alexander-Whitney diagonal approximation as a (graded) derivation to obtain a monoid $\mathbb{A}(K)$ in the monoidal category of dg $C_0(K)$-bicomodules. 
\end{enumerate}

We consider two free resolutions of $\mathbb{A}(K)$ as a differential graded $\mathbb{A}(K)$-bimodule. On one hand, we have the \textit{bar resolution}, which has as underlying graded $\mathbb{A}(K)$-bimodule
\[\mathcal{B}(\mathbb{A}(K))= \mathbb{A}(K) \square \big(\bigoplus_{n=0}^{\infty} s^{+1}\mathbb{A}(K)\overbrace{\square \cdots \square}^{n} s^{+1}\mathbb{A}(K)\big) \square \mathbb{A}(K),\]
where $s^{+1}$ denotes the functor that shifts the grading up by $1$. The differential of $\mathcal{B}(\mathbb{A}(K))$ is defined in terms of the internal differential of $\mathbb{A}(K)$ and the concatenation product. On the other, there is a smaller resolution with underlying graded $\mathbb{A}(K)$-bimodule given by
\[ \mathcal{Q}( \mathbb{A}(K)) = \mathbb{A}(K) \square C(K) \square \mathbb{A}(K),\]
where $C(K)$ is the free graded $R$-module generated by simplices in $K$. The differential of $\mathcal{Q}( \mathbb{A}(K))$ is defined in terms of the simplicial boundary and Alexander-Whitney diagonal approximation of $C(K)$, the natural degree $-1$ inclusion map $C(K) \to \mathbb{A}(K)$, and the product and the internal differential of $\mathbb{A}(K)$. 

The following statement gives a combinatorial interpretation of each of the algebraic constructions described above. It uses the notion of a $\mathbf{P}$-\textit{complex} for a set of polytopes $\mathbf{P}$ (Definition \ref{P-complex}), which is a version of a polytopal complex but allowing a cell to have face identifications, generalizing the notion of $\mathbf{\Delta}$-complex \cite[Appendix]{hatcher2002algebraic}. 

\begin{customthm}{2}\label{theorem2} Let $K$ be a simplicial complex with a coherent ordering of the vertices of each simplex. 
\begin{enumerate}
    \item $\mathbb{A}(K)$ is isomorphic to the chains on a $\mathbf{I}$-complex $\mathcal{A}_K$, where $\mathbf{I}=\{I^k\}_{k \in \mathbb{N}_0}$ and $I^k$ denotes the $k$-dimensional cube.
    \item $\mathcal{Q}(\mathbb{A}(K))$ is isomorphic to the chains on a $(\mathbf{F}\times \mathbf{I})$-complex $\mathcal{Q}_K$, where  $\mathbf{F}\times \mathbf{I} = \{F_p \times I^q\}_{p,q \in \mathbb{N}_0}$, with $F_p$ being the $p$-dimensional \emph{freehedron}.
    \item $\mathcal{B}(\mathbb{A}(K))$ is isomorphic to the chains $C(\mathcal{B}_K)$ on a $(\mathbf{\Delta}\times\mathbf{I})$-complex $\mathcal{B}_K$,  where $\mathbf{\Delta}\times\mathbf{I}= \{\Delta^p \times I^q\}_{p,q\in \mathbb{N}_0}$. 
  \item   There is an injective quasi-isomorphism of $\mathbb{A}(K)$-bimodules
\[\alpha\colon\mathcal{Q}(\mathbb{A}(K))\to\mathcal{B}(\mathbb{A}(K)) \cong C(\mathcal{B_K})\] such that, for any simplex $\sigma \in K_n$, \[\alpha(\operatorname{id}_{s(\sigma)} \square \sigma \square \operatorname{id}_{t(\sigma)})\in C(\mathcal{B_K})\] is a sum of prisms whose union forms a polytope combinatorially equivalent to $G_n$. 

\end{enumerate}
 \end{customthm}
By \textit{chains} in each of the above statements, we mean the $R$-chain complex generated by the cells in each respective $\mathbf{P}$-complex with boundary differential determined by an appropriate orientation of each cell.   

Statement (1) goes back to the work of Baues and Rivera-Zeinalian \cite{Bau80, rivera2018cubical}. This observation underlies the \textit{Baues generalized problem} as posed in \cite{BKS}. The key observation is that the nonempty face poset of $I^{k-1}$ is isomorphic to the poset of necklaces in $\Delta^k$ joining the first and last vertices.

The \textit{freehedra} (also known in the literature as \textit{Hochschild polytopes}) form a family of polytopes appearing when constructing combinatorial models for free loop spaces as described in \cite{S09, rivera-saneblidze}. Statement (2) essentially follows from this work. The key observation is that the nonempty face poset of $I^{p-1} \times F_k \times I^{q-1}$ is isomorphic to the poset of \textit{marked necklaces} inside $\Delta^p \times \Delta^k \times \Delta^q$. The combinatorial and algebraic meaning of freehedra has been further studied in \cite{C25, Com21,  poliakova2020cellular, PP25}. 

Statement (3) is a consequence of (1) and the simplicial nature of the bar resolution. 

Goodwillie polytopes finally appear in statement (4) as the union of certain prisms in the cell complex realizing bar resolution. The quasi-isomorphism $\alpha$ may be regarded as a manifestation of \textit{Koszul duality} and this statement can be understood as providing a combinatorial picture for it. The authors discovered Goodwillie polytopes through a similar manifestation of this fact in the context of singular chains on topological spaces motivated by constructing algebraic models for loop spaces in \cite{rivera_cyclic_2024}.

We end this introductory discussion by posing three questions.
\subsection*{Questions}
\begin{enumerate}
\item Can the face numbers of prismatic polytopes be characterized? \label{Q1}

 \item Can all families of polytopes satisfying the properties (1)-(4) in Theorem \ref{theorem1}, or all self-dual prismatic polytopes at a fixed dimension, be classified? \label{Q2}

\item What is the deeper meaning of the self-duality of the Goodwillie polytopes in the light of Theorem \ref{theorem2}? Is the self-duality a consequence of applying a construction that preserves such property to the family of simplices? Can this construction be iterated to construct new families of polytopes with interesting properties? \label{Q3}
\end{enumerate}

\subsection*{Creative process, AI usage, and acknowledgments} 
The family of Goodwillie polytopes was discovered during DT's PhD thesis work, which constructed a model for the $S^1$-equivariant homology of the free loop space of a (possibly non-simply connected) topological space in terms of its coalgebra of singular chains; the core of this work appears in \cite{rivera_cyclic_2024}. The properties of Goodwillie polytopes studied in the present article are not used in achieving the main goal of that previous work. Rather, our interest in the combinatorics and deeper meaning of these polytopes grew out of long, slow, AI-free, tangential discussions, together with sporadic interactions with combinatorialists, aimed at understanding Koszul duality in terms of pictures.
 
A few weeks after \cite{rivera_cyclic_2024} was posted, Frédéric Chapoton shared with us a SageMath implementation computing the face vector of the Goodwillie polytopes (and suggested a correction to the original formula for the coordinates of their vertices). These first computations suggested that the face vector of each Goodwillie polytope is palindromic. After analyzing some low dimensional examples, this led us to conjecture the self-duality proven in the present article. We thank F. Chapoton for the fruitful interactions that sparked our interest in the combinatorial structure of the Goodwillie polytopes.
 
Verifying the self-duality conjecture involved computational experiments. For this, DT used ChatGPT 5.6 Sol on Codex to generate code for the necessary computations, experiments, and visualizations, some of which are available in the following interactive website: \href{https://goodwillie-polytopes-visualizations.danieltolosa.workers.dev}{A hands-on guide to Goodwillie polytopes} \cite{tolosa_goodwillie_visualizations}. AI tools were also used to search the literature; in particular, they pointed to useful references on regular polytopal subdivisions, providing the right language to state a decomposition of Goodwillie polytopes into prisms. In all cases, the outputs of these tools were critically examined, digested, and discussed by the authors. All the writing was produced by the authors. 

The context in which the in-person discussions, writing, and thinking behind this project took place also deserves comment, as it shaped the creative process in an important way. A significant part of this project was carried out in different cities in Colombia (Cali and Bogotá). DT gave a talk on these results at the \textit{Second Colombian Geometry and Topology Meeting} (\href{https://www.ecogyt.com}{ECOGyT}), held in Cali, Colombia, from August 3--14, 2026, of which MR was one of the organizers. On August 10, a major 7.4-magnitude earthquake struck Cali and other regions in the country. Somehow the organizers and participants decided to reconvene the event in a different format and location.
 
It was strange to hold mathematical discussions amid so much human suffering. And yet many of us witnessed how mathematical interaction and the pursuit of human understanding can offer grounding and direction in moments of need and confusion. We dedicate this work to the victims and survivors of that catastrophe. 

The interactions of those days had a particular weight and pace that is hard to put into words, but that shaped the perspective of this work. In particular, we would like to thank Camilo Arias-Abad and Federico Ardila for stimulating questions and discussions during the event. 

Two weeks after the earthquake, MR presented this work in a celebration of 70 years of the \textit{Universidad Nacional de Colombia (UNAL)} Mathematics Department in Bogotá, in a talk accessible to a wide audience. Fruitful discussions with logicians, geometers, combinatorialists, and mathematical physicists followed. MR would like to thank the organizers of such a stimulating event and the \textit{UNAL}, Bogotá for its excellent working conditions and hosting a sabbatical visit, during which most of this work was written and polished.
 
It is difficult to envision what shape and direction this project would have taken had the initial research in \cite{rivera_cyclic_2024} been conducted with the aid of AI tools and guided primarily by the goal of proving its main result, or without the sustained in-person mathematical discussions, talks, and seminars in Colombia that accompanied it. Much of the understanding and insight that eventually led to this work emerged from that earlier, exploratory phase. At the same time, AI tools proved useful in the later stages of the project, particularly for the computational tasks involved in testing and verifying our conjectures.

We include these reflections because we believe that the process through which mathematical understanding develops deserves some attention, beyond the mathematical results that eventually emerge from it.

Finally, MR acknowledges support by NSF Grant DMS 2405405 and by Shota Rustaveli National Science Foundation of Georgia (SRNSFG) [grant number FR-23-5538]. 

\subsection*{Organization}

The paper is organized as follows. In Section~\ref{preliminaries}, we review the necessary preliminaries from polytopal combinatorics. In Section~\ref{sec:goodwillie}, we introduce the Goodwillie polytopes. In Section~\ref{sec:regular-subdivision}, we describe their decomposition into prisms and prove part~(2) of Theorem~\ref{theorem1}. In Section~\ref{sec:self-duality}, we establish parts~(1) and~(3) of Theorem~\ref{theorem1}. In Section~\ref{sec:face-enumeration}, we derive a closed formula for the face numbers and the $f$-polynomial of $G_n$, thereby proving part~(4) of Theorem~\ref{theorem1}. Finally, in Section~\ref{sec:resolutions}, we define $\mathbb{A}(K)$, construct the two resolutions, and prove Theorem~\ref{theorem2}.

\section{Preliminaries} \label{preliminaries}

\subsection{Polytopes, face lattices, and subdivisions}

We briefly recall the terminology from polytopal combinatorics used
throughout the paper following
\cite{ziegler2012lectures} and
\cite[Chapter~2]{de_loera_triangulations_2010}.

\begin{definition}
    A \textit{convex polytope} is the convex hull
    \[
        P=\conv(A)
    \]
    of a finite point configuration $A\subseteq\mathbb R^n$. The \textit{dimension of} $P$ is the dimension of its affine span in $\mathbb{R}^n$.
\end{definition}

\begin{definition}
    A \textit{face} of a polytope $P \subseteq \mathbb{R}^n$ is a subset of the form
    \[P \cap \{(x_1,\ldots,x_n) \in \mathbb{R}^n : c_1x_1 + \cdots + c_nx_n = C\}\]
    for some $(c_1,\ldots,c_n) \in \mathbb{R}^n$ and $C \in \mathbb{R}$ such that $c_1x_1 + \cdots +c_nx_n \leq C$ is satisfied for all $(x_1,\ldots, x_n) \in P$. The \textit{dimension of a face} is defined as the dimension of its affine hull. Faces of dimension $0,1,$ and $\dim P-1$ are called \textit{vertices, edges}, and \textit{facets}, respectively. By convention, the empty face has dimension $-1$. The \textit{face lattice} $\mathcal F(P)$ is the poset of all faces of
    $P$, ordered by inclusion.
\end{definition}
There is a convenient equivalent way of describing polytopes and their faces via half-spaces and supporting hyperplanes.
\begin{theorem}[Minkowski--Weyl theorem for polytopes]\label{thm:minkowsky--weyl}
    Let $P\subseteq \R^n$ be non-empty. The following are equivalent:
    \begin{enumerate}
        \item There exists a finite set $A\subseteq \R^n$ such that
        \[
            P=\conv(A).
        \]
        \item There exist finitely many affine maps
        $\ell_1,\dots,\ell_m\colon \R^n\to\R$ such that
        \[
            P=\bigcap_{i=1}^m
            \left\{x\in\R^n\ \middle|\ \ell_i(x)\leq 0\right\},
        \]
        and this intersection is bounded.
    \end{enumerate}
    Thus, a polytope may equivalently be defined as the convex hull of
    finitely many points or as a bounded intersection of finitely many
    closed affine half-spaces.
\end{theorem}

A hyperplane $H\subseteq\R^n$ is said to \emph{support} a polytope
$P$ if $P\cap H$ is a non-empty face of $P$.

\begin{definition}
    Two polytopes $P$ and $Q$ are said to be \textit{combinatorially equivalent} if there is an isomorphism of posets $\mathcal F(P)\cong\mathcal F(Q)$ between their face lattices. The combinatorial type of $P$ is its equivalence class under this relation.
\end{definition}

Thus two combinatorially equivalent polytopes need not be affinely or
geometrically equivalent. They need only have the same face-containment structure.

\begin{definition}

    For a polytope $Q$, its \textit{relative interior} $\operatorname{relint}(Q)$
is its interior taken inside its affine hull
$\operatorname{aff}(Q)$.    
\end{definition}

\begin{definition}\label{P-complex}
Let $\mathbf{P}$ be a set of convex polytopes, closed under taking faces up to affine
isomorphism. A \emph{$\mathbf{P}$-complex} is a CW complex $X$ together with, for every cell
$e$ of $X$, a polytope $P_e\in\mathbf{P}$ and a characteristic map \[\chi_e\colon P_e\to X\] such that:
\begin{enumerate}
  \item the restriction of $\chi_e$ to the relative interior of $P_e$ is a homeomorphism
        \[\chi_e|_{\operatorname{relint}(P_e)}\colon\operatorname{relint}(P_e)\xrightarrow{\ \cong\ }e;\]
  \item for every nonempty proper face $F$ of $P_e$, there are a cell $e'$ of $X$ and an affine
        isomorphism $\phi\colon F\to P_{e'}$ with \[\chi_e|_F=\chi_{e'}\circ\phi.\]
\end{enumerate}
\end{definition}

When $\mathbf{P}=\mathbf{\Delta}=\{\Delta^k\}_{k\in \mathbb{N}_0}$, the above definition agrees with that of a $\mathbf{\Delta}$-complex in \cite[Appendix]{hatcher2002algebraic}.

\begin{definition}
    A finite collection $\Sigma$ of polytopes in a common Euclidean
    space is a \emph{polytopal complex} if:
    \begin{enumerate}
        \item every face of every $Q\in\Sigma$ also belongs to
        $\Sigma$;
        \item for every $Q,Q'\in\Sigma$, the intersection
        $Q\cap Q'$ is a face of both $Q$ and $Q'$, possibly the
        empty face.
    \end{enumerate}
\end{definition}

\begin{definition}\label{def:polytopal-subdivision}
    Let $P$ be a polytope. A \emph{polytopal subdivision} of $P$ is a
    polytopal complex $\Sigma$ such that
    \[
        \bigcup\limits_{Q\in \Sigma} Q=P
    \]
    and every maximal cell of $\Sigma$ has dimension $\dim P$.\\
    If, additionally, every nonempty cell $Q$ is combinatorially equivalent to a prism, i.e.~a Cartesian product $\Delta_p \times I^q$, for some $p,q \in \mathbb{N}_0$ possibly varying from cell to cell, we say that $\Sigma$ is a \textit{prismatic subdivision} of $P$. 
\end{definition}
Equivalently, one may specify only the maximal cells of a subdivision,
provided that all of their faces are understood to be included. We will use this description often.

We now recall a standard construction that produces polytopal
subdivisions (\ref{def:polytopal-subdivision}). Let $A\subseteq\mathbb R^n$ be a finite point configuration affinely spanning $\mathbb R^n$, and let
\[
    h:A\longrightarrow\mathbb R
\]
be any function. Define the \textit{lifted point configuration} by
\[
    A^{h}
    :=
    \{(a,h(a)):a\in A\}
    \subseteq\mathbb R^{n+1},
\]
whose convex hull is the \textit{lifted polytope}
\[
    P^{h}
    :=
    \conv(A^{h}).
\]
We call $h$ a \textit{height} function for $A$. 

We shall denote by
\[
    \rho:\mathbb R^{n+1}\longrightarrow\mathbb R^n,
    \qquad
    \rho(x,z)=x,
\]
the projection that forgets the last coordinate.

\begin{definition}\label{def:lower-face}
    Let $A$ be a finite point configuration affinely spanning $\R^n$. A face $F$ of the lifted polytope $P^{h}$ is a \textit{lower face} if
    there exists an affine function
    \[
        \ell:\mathbb R^n\longrightarrow\mathbb R
    \]
    such that
    \[
        h(a)\geq\ell(a)
        \qquad\text{for every }a\in A,
    \]
    and
    \[
        F
        =
        P^{h}
        \cap
        \{(x,\ell(x))\in\mathbb R^{n+1}\}.
    \]
\end{definition}
We will use the following standard result
\cite[Lemmas~2.3.4 and 2.3.11]{de_loera_triangulations_2010}.
\begin{theorem}
    \label{thm:regular-subdivision}
    Let $A\subseteq\mathbb R^n$ be a finite point configuration whose affine span is
    $\mathbb R^n$, let $h:A\to\mathbb R$ be a function, and denote $P=\conv(A)$. Then the collection
    \[
        \Sigma_h
        :=
        \left\{
            \rho(F):F \text{ is a lower face of } P^{h}
        \right\}
    \]
    is a polytopal subdivision of $P$. In particular,
    \begin{enumerate}
        \item the cells of $\Sigma_h$ cover $P$, namely
        \[
            \bigcup_{Q\in\Sigma_h}Q=P \ ,
        \]
        \item the intersection of any two cells is a (possibly empty) common face of
        both,
        \item distinct maximal cells have disjoint relative interiors, and
        \item if $\dim(P^h)=n+1$ then, for a lower face $F$ of $P^h$, $F$ is a lower facet of $P^h$ if and only if $\rho(F)$ is a maximal cell of $\Sigma_h$.
    \end{enumerate}
\end{theorem}
A polytopal subdivision arising from the construction we just described is said to be \textit{regular}.

\section{\texorpdfstring{The Goodwillie polytopes $\G$}{The Goodwillie polytopes G n epsilon}}\label{sec:goodwillie}
Goodwillie polytopes are a family $\mathcal{G}= \{G_0, G_1, G_2, G_3, \dots\}$, consisting of one convex polytope in each dimension.
$G_0$ and $G_1$ are the combinatorially equivalent to a point and the interval $[0,1] \subset \R$, respectively; $G_2$ is the combinatorially equivalent to the convex hull of the points $(0,0),(0,1),(1,0),(1,1),(\frac{1}{2},-1) \in \R^2$, i.e.~ a pentagon. For these low dimensions, the statements regarding a regular prismatic subdivision, self-duality, and $f$-vectors, are immediate. We will therefore consider the case $n\geq 3$ throughout. 

We will define geometric convex polytopes $\G$ that depend on a parameter $\epsilon >0$ and show that the combinatorial type of $\G$ is stable in the range $0< \epsilon < \min (\frac{1}{2}, \frac{1}{n-2})$. Once the stable range for $\epsilon$ has been established (Corollary \ref{cor:stable-range}) we will drop the $\epsilon$ from the notation and denote by $G_n$  both the geometric and combinatorial polytopes.

Each $\G$ arises as the convex hull of the union of the vertices of the $n$-cube and an additional set of points as defined next.

\begin{definition}\label{def:G_n}
    Let $n\geq 3$ and suppose $0<\epsilon< \min \set{\frac{1}{2},\frac{1}{n-2}}$.
    For $t\in\{0,1\}$ and $A\subseteq \set{2,...,n}$, define $c_{t,A}\in\R^n$ by
\[
    (c_{t,A})_1=t,
    \qquad
    (c_{t,A})_i=
    \begin{cases}
        1,& i\in A,\\
        0,& i\notin A,
    \end{cases}
    \quad i\in \set{2,...,n}.
\]
Thus the points $c_{t,A}$ are precisely the vertices of the cube
$[0,1]^n$, with $A$ recording which coordinates other than the first
are equal to $1$. For $r\in \set{2,...,n}$ and $B\subseteq \set{2,...,n}\setminus\{r\}$, define
$u_{r,B}^{\epsilon}\in\R^n$ by
\[
    (u_{r,B}^{\epsilon})_1=\frac{r-1}{n},
    \qquad
    (u_{r,B}^{\epsilon})_r=-1,
\]
and
\[
    (u_{r,B}^{\epsilon})_i=
    \begin{cases}
        1-\epsilon,& i\in B,\\
        \epsilon,& i\notin B,
    \end{cases}
    \quad i\in \set{2,...,n}\setminus\{r\}.
\]
Here, the set $B$ is recording which coordinates are equal to $1-\epsilon$. We say $u^\epsilon_{r,B}$ \textit{is of type $r$}, to indicate that its $r$-th coordinate is the distinguished coordinate with entry $-1$.
We call the $c_{t,A}$ and $u^\epsilon_{r,B}$ \textit{cube points} and \textit{added points}, respectively, and refer to their union as the \textit{defining points} of $\G$.

    The $(n,\epsilon)$-Goodwillie polytope is defined by
\[
\begin{aligned}
    G_n^\epsilon
    &\coloneqq \conv\bigl(
        \{c_{t,A}:t\in\{0,1\},\ A\subseteq \set{2,...,n}\}\\
    &\qquad{}\cup
        \{u_{r,B}^{\epsilon}:r\in \set{2,...,n},\ B\subseteq \set{2,...,n}\setminus\{r\}\}
    \bigr).
\end{aligned}
\]
\end{definition}

\begin{figure}[H]
    \centering
    \includegraphics[width=0.5\textwidth]{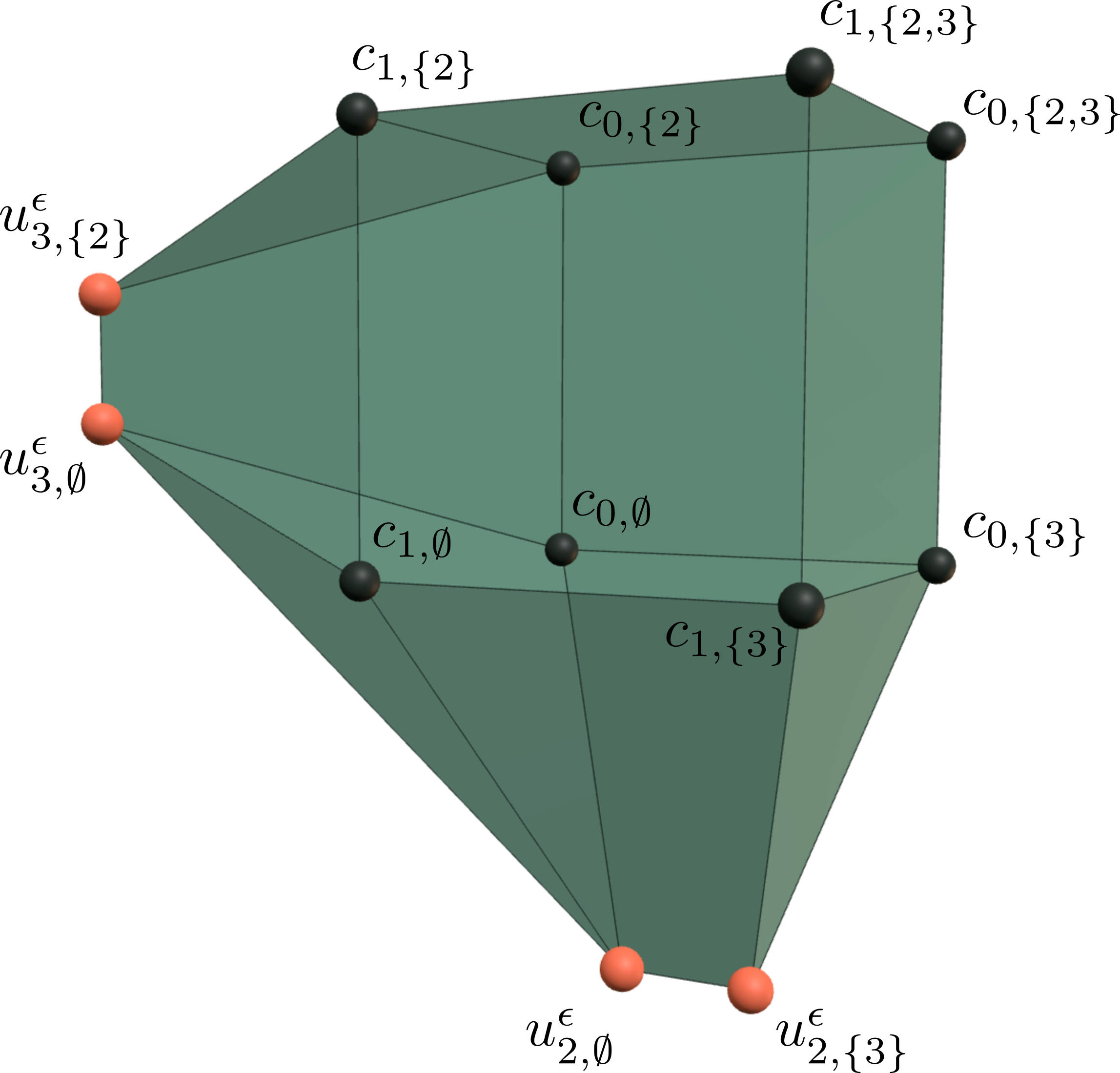}
\caption{The Goodwillie polytope $G_3^\epsilon$ for
$\epsilon=\frac13$, with cube vertices in black and added
vertices in orange. As $\epsilon$ increases towards $\frac12$,
the cube vertices remain fixed, while the two added vertices
of each type approach one another: their differing coordinate
values are $\epsilon$ and $1-\epsilon$, so their distance is
$1-2\epsilon$. Each pair coincides at $\epsilon=\frac12$. If $\epsilon = 0$, the face structure degenerates as distinct faces become coplanar, e.g.~the face $\conv\left(c_{1,\set{3}}, c_{0,\set{3}}, u_{2,\set{3}}^\epsilon\right)$ and the cube face $\conv\left(c_{1,\set{2,3}}, c_{0,\set{2,3}}, c_{0,\set{3}}, c_{1,\set{3}}\right)$ merge into a single pentagonal face.}
    \label{fig:G_3}
\end{figure}

The reason for requiring $0< \epsilon < \frac{1}{2}$ is clear from the definition, as equality in either side would change the combinatorial type of $\G$ (see Fig \ref{fig:G_3} for a visual example).
\begin{remark}
    Goodwillie polytopes admit a realization as integral polytopes (convex polytopes whose vertices have integer coordinates) for example by setting $\epsilon = \frac{1}{n}$ and scaling $\G$ by $n$.
\end{remark}

\begin{proposition}
$\G \subset \R^n$ is an $n$-dimensional convex polytope whenever $0<\epsilon< \min \set{\frac{1}{2},\frac{1}{n-2}}$.
\end{proposition}
\begin{proof}
    The fact that it is a polytope is obvious as $\G$ was defined as the convex hull of a configuration of points in $\R^n$. The dimension is also clear, as the points of the $n$-cube, which are a subset of the point configuration defining $\G$, already span $\R^n$ affinely, and this is the largest affine dimension that a configuration of points in $\R^n$ can attain.
\end{proof}
\Needspace{7\baselineskip}
\subsection{\texorpdfstring{The prisms $\P$}{The prisms P S epsilon}}
\begin{definition}\label{def:P}
    For $S\subseteq \set{2,\dots, n}$, define $\P \subset \G$ to be the convex hull of the points
    \begin{enumerate}
        \item $c_{t,A}$ such that $A \cap S = \emptyset$,
        \item $u^\epsilon_{r,B}$ such that $r \in S$ and $B\cap S = \emptyset$.
    \end{enumerate}
\end{definition}

\noindent These conditions force every coordinate in $S$ to be \textit{low}, i.e.~ to be $0$ for cube vertices, and $\epsilon$ or $-1$ for added vertices. We therefore think of $S$ as the \textit{coordinates that are forced to be low} for the points in the configuration defining $\P$. 

Label the elements of $S \subseteq \set{2,\dots,n}$ so that $S = \set{ i_1 < i_2< \cdots < i_k}$ and consider the projection to the coordinates of $\R^n$ indexed by $S\cup \set{1}$
\begin{equation}\label{eqn:proj}
\pi_{\Delta}: \P \rightarrow \R^{k+1} \ , \quad \pi_{\Delta}(x_1 , \dots , x_n) = (x_1,x_{i_1}, \dots, x_{i_k}).
\end{equation}
The cube points $c_{t,A}$ map to either
\[
\zeta_0 = (0,0,\dots,0) \text{ or } \zeta_1 = (1,0,\dots,0) ,
\]
and, for each $r \in S$, every added point of type $r$ in $\P$ maps to
\[
\upsilon_r = \left(\frac{r-1}{n}, \epsilon, \dots, \epsilon, -1, \epsilon, \dots, \epsilon \right)
\]
where the $-1$ is in the coordinate obtained from projecting the $r$-th coordinate of $\R^n$ under $\pi_{\Delta}$.

\begin{proposition}\label{prop:proj_is_simplex}
    If $0 < \epsilon < \min\set{\frac{1}{2},\frac{1}{n-2}}$, then the set of $k+2$ points \[ \{\zeta_0, \zeta_1,\} \cup \{ \upsilon_r\ \vert \  r \in S\} \]
    is affinely independent. Consequently, $\conv (\{\zeta_0, \zeta_1,\} \cup \{ \upsilon_r\ \vert \  r \in S\})$ is a $k+1$-simplex.
\end{proposition}
\begin{proof}
    The result holds trivially for $k\in \set{0,1}$. Assume $k\geq 2$. A set of points $\set{p_1,\dots , p_t}$ in a vector space $V$ is affinely independent if and only if the points $\set{p_2 - p_1, p_3-p_1,\dots , p_t-p_1}$ are linearly independent in $V$. Hence, it suffices to check that the elements of $\{\zeta_1\} \cup \{ \upsilon_r\ \vert \  r \in S\}$ are linearly independent in $\R^{k+1}$. This is equivalent to finding the non-vanishing conditions for the determinant
    \[
    \begin{vmatrix}
       \ \ 1&\frac{{{i_1}}-1} {n}&\frac{{{i_2}}-1}{n}&\cdots&\frac{{{i_k}}-1}{n}\ \\
       \ \ 0&-1&\epsilon& \cdots&\epsilon\ \\
       \ \ 0&\epsilon&-1&\cdots&\epsilon\ \\
       \ \ \vdots&\vdots & \vdots& \ddots & \vdots\ \\
        \ \ 0&\epsilon&\epsilon&\cdots&-1\ 
    \end{vmatrix} =
        \begin{vmatrix}
       \ -1&\epsilon& \cdots&\epsilon\ \\
       \ \epsilon&-1&\cdots&\epsilon\ \\
       \ \vdots&\vdots & \ddots & \vdots\ \\
        \ \epsilon&\epsilon&\cdots&-1\ 
    \end{vmatrix} = (-1)^k(1+\epsilon)^{k-1}(1-(k-1)\epsilon).
    \]

Thus, for $k\geq2$, the determinant is zero if and only if $\epsilon =\frac{1}{k-1} $. %
The smallest value of $\epsilon$ for which this happens, for a fixed $n$, is when $|S|=k=n-1$ and $\epsilon = \frac{1}{n-2}$. Therefore, if $0<\epsilon<\frac{1}{n-2}$, the result holds. 
\end{proof}
We now describe the fibers of $\pi_\Delta$. A point $q \in  \conv (\{\zeta_0, \zeta_1,\} \cup \{ \upsilon_r\ \vert \  r \in S\})$ has barycentric coordinates
\[q= z_0 \zeta_0 + z_1 \zeta_1 + \sum\limits_{t=1 }^{k} w_{i_t} \upsilon_{i_t}, \] 
where $z_0,z_1,w_{i_1}, \dots, w_{i_k} \geq 0$ and $z_0 + z_1 + w_{i_1} + \cdots + w_{i_k} = 1$.\\
Define the affine function
\[
\lambda_S (q) = \sum_{r \in S} w_r.
\]
This function satisfies $\lambda_S(\zeta_0)=\lambda_S(\zeta_1)=0$, $\lambda_S(\upsilon_r)=1$, and $0\leq \lambda_S(q) \leq 1$.
Write $\set{2,...,n}\setminus S = \set{j_1,\dots,j_{n-k-1}}$. All points $x \in \P$ satisfy
\[
\begin{gathered}
\epsilon \lambda_S(\pi_{\Delta}(x)) \leq x_j \leq 1-\epsilon \lambda_S(\pi_{\Delta}(x)),\\
\forall j\in \set{2,...,n}\setminus S.
\end{gathered}
\] 
Therefore, for any $q \in \pi_\Delta(\P)$, the fiber $\pi_\Delta^{-1}(q)$ is contained in a cube $[\epsilon \lambda_S(q) , 1-\epsilon \lambda_S(q)]^{n-k-1}$. Conversely, fix \(q\) with barycentric coordinates \(z_0,z_1,(w_r)_{r\in S}\). For each \(A\subseteq \set{2,...,n}\setminus S\), form the point
\[
x=z_0c_{0,A}+z_1c_{1,A}
+\sum_{r\in S}w_r u^\epsilon_{r,A}\in P_S^\epsilon.
\]Then \(\pi_\Delta(x)=q\), and for every \(j\in \set{2,...,n}\setminus S\),
\[
x_j=
\begin{cases}
1-\epsilon\lambda_S(q),&j\in A,\\
\epsilon\lambda_S(q),&j\notin A.
\end{cases}
\]
Thus the fiber contains every corner of the claimed cube, and hence the entire cube by convexity. This cube is combinatorially equivalent to an $(n-k-1)$-cube, as long as $\epsilon < \frac{1}{2}$.\\
We now give an equivalent description of $\P$ as a subset of $\fancy{H} \times \R^{n-k-1}$, where barycentric coordinates identify $\pi_\Delta(\P)$ with a subset of
\[
 \fancy{H} = \set{ (z_0,z_1,w_{i_1},\dots,w_{i_k}) \in \R^{k+2} \ \vert \ z_0+z_1+w_{i_1}+\cdots+ w_{i_k} = 1},
\]
which is an affine hyperplane of $\R^{k+2}$ of dimension $k+1$.\\ 
Write $q = (z_0,z_1,w_{i_1},\dots,w_{i_k}) \in \fancy{H}$, and $y = (y_{j_1}, y_{j_2}, \dots, y_{j_{n-k-1}})$, and let $\fancy{P}^\epsilon_{S}$ be the subset of $\fancy{H} \times \R^{n-k-1}$ of points $(q,y)$ such that
\begin{enumerate}
    \item $ z_0 \geq 0$,
    \item $z_1 \geq 0$,
    \item $w_r \geq 0$ \text{for each $r \in S$}
    \item $y_j - \epsilon \lambda_S(q) \geq 0$ \text{for all $j \in \set{2,...,n} \setminus S$,} \text{ and }
    \item $1- y_j - \epsilon \lambda_S(q) \geq 0$ \text{for all $j \in \set{2,...,n} \setminus S$}. 
\end{enumerate}

\begin{proposition}\label{prop:prism_facets}
    The subset $\fancy{P}^\epsilon_S \subset \fancy{H} \times \R^{n-k-1}$ defines an $n$-dimensional convex polytope mapped affinely isomorphically to $\P$ under the natural map $\fancy{H}\times \R^{n-k-1} \xrightarrow{\cong} \R^{n}$. The vertices of $\P = \fancy{P}^\epsilon_S$ are the elements of the configuration of points defining $\P$ (\ref{def:P}).
\end{proposition}
\begin{proof}
    After identifying $\fancy{H}$ with the affine hull of $\pi_\Delta(\P)$, and $\R^{n-k-1}$ with the coordinates of $R^n$ indexed by $\set{2,...,n}\setminus S$, it is clear that $\P = \fancy{P}_S^\epsilon$. The latter is defined by a bounded non-empty intersection of half-spaces that contains a cartesian product of a $(k+1)$-dimensional simplex and a cube $[\epsilon,1-\epsilon]^{n-k-1}$, with $0<\epsilon<\frac{1}{2}$, and is therefore an $n$-dimensional convex polytope (\ref{thm:minkowsky--weyl}). Lastly, let \((q,y)\in\fancy{P}_S^\epsilon\)
    be a vertex. If some coordinate $y_j$ lies strictly between $\epsilon\lambda_S(q)$ and $1-\epsilon\lambda_S(q)$, we can increase or decrease it by the same sufficiently small amount, keeping $q$ and all other coordinates fixed. Both resulting points lie in the prism, and $(q,y)$ is their midpoint, contradicting that it is a vertex. Thus every $y_j$ must be an endpoint of its interval. Similarly, one can show that \(q\) must be a vertex of $\pi_\Delta(\P)$. Conversely, suppose that \(q\) is a vertex of the barycentric simplex and every \(y_j\) is an endpoint of its interval. If \((q,y)\) is the midpoint of two points of \(\fancy{P}_S^\epsilon\), then their projections must both equal \(q\), since \(q\) is a vertex of the simplex. Both points therefore belong to the fiber over \(q\). But \(y\) is a vertex of this cubical fiber, so the two points must both equal \((q,y)\). Hence \((q,y)\) is a vertex. This gives precisely
    \[   \{c_{t,A} \ \vert \  A \cap S = \emptyset\} \cup \{ u^\epsilon_{r,B} \ \vert \ r \in S \ \text{and} \ B\cap S = \emptyset\}
    \]
    with $c_{t,A} , u^\epsilon_{r,B}$ as defined in \ref{def:G_n}.
\end{proof}
\begin{corollary}\label{cor:facets_P}
    The complete list of facets of $\P$ is
\begin{align}
 &\mathbf{B}_{\zeta_0,S}= \{ (q,y) \in \P \ | \  z_0 = 0\} ,\\
&\mathbf{B}_{\zeta_1,S}= \{ (q,y) \in \P\  | \ z_1=0\} ,\\
&\mathbf{B}_{\upsilon_r,S}= \{ (q,y) \in \P\  |\  w_r =0\}, \quad  r\in S\ ,\\
&\mathbf{L}_{j,S}= \{ (q,y)\in \P \  |  \ y_j - \epsilon \lambda_S(q) = 0 \}, \quad j \in \set{2,...,n}\setminus S,\\ 
&\mathbf{U}_{j,S}=\{ (q,y) \in \P  \ | \  1- y_j - \epsilon \lambda_S(q)=0\},\quad  j \in \set{2,...,n}\setminus S .
\end{align}
   Consequently, $\P$ has $2n-k$ facets.
\end{corollary}
We will refer to these families as $\mathbf{B}$-, $\mathbf{L}$-, and $\mathbf{U}$-facets.
\begin{theorem}\label{thm:prism-comb-type}
    For $n \geq 3, \  0< \epsilon < \min\{\frac{1}{2},\frac{1}{n-2}\}$, and $S \subseteq \set{2,\dots,n}$, $\P$ is combinatorially equivalent to
    \[
    \Delta^{\abs{S}+1} \times I^{n-1-\abs{S}}.
    \]
    
\end{theorem}
    \begin{proof}
    Let $k=\abs{S}$ and $m=n-1-\abs{S}$. We begin by describing the face lattice of
    \[
        \Delta^{k+1}\times I^m.
    \]
    A nonempty face of a Cartesian product is uniquely determined by a choice of a nonempty face of each factor. In our case, this amounts to choosing a nonempty face $F$ of the simplex $\Delta^{k+1}$, together with a face of the cube $I^m$.\\
    In barycentric coordinates, a face of a simplex is determined by specifying which coordinates are constrained to be zero: setting the coordinate associated to a vertex $v$ to $0$ amounts precisely to restricting to the face opposite $v$. On the other hand, a face of the cube
    \[
        I^m=\prod_{j\in \set{2,...,n}\setminus S}I
    \]
    is determined by choosing, in each interval factor, the lower endpoint $\{0\}$, the upper endpoint $\{1\}$, or the whole interval $I$. Therefore, the nonempty faces of $\Delta^{k+1}\times I^m$ are labeled exhaustively and uniquely by tuples
    \[
        \left(D,(E_j)_{j\in \set{2,...,n}\setminus S}\right),
    \]
    where $D$ is a nonempty face of $\Delta^{k+1}$ and $E_j\in\bigl\{\{0\},\{1\},I\bigr\}$. These are the elements of the face lattice. Their partial order is given by coordinatewise containment:
    \[
        \left(D,(E_j)\right)\subseteq \left(D',(E'_j)\right)
    \]
    if and only if
    \[
        D\subseteq D'
        \qquad\text{and}\qquad
        E_j\subseteq E'_j
        \quad\text{for every }j\in \set{2,...,n}\setminus S.
    \]
    
    We will show that the nonempty faces of $\P$ admit exactly the same description.
    Any face of a polytope is the intersection of all the facets containing it. Thus, if $F$ is a nonempty face of $\P$, then
    \begin{equation}\label{eqn:face_facets}
        F
        =
        \bigcap_{\substack{H\in\operatorname{Facets}(\P)\\F\subseteq H}}H,
    \end{equation}
    with the empty intersection understood as $\P$.
    By \ref{prop:proj_is_simplex}, the projection $\pi_\Delta(\P)$ is combinatorially equivalent to $\Delta^{k+1}$. By \ref{cor:facets_P}, the facets $\mathbf B_{\bullet,S}$ of $\P$ are in correspondence with the facets of $\pi_\Delta(\P)$. Hence, the collection of $\mathbf{B}$-facets containing $F$ determines a unique nonempty face $\pi_\Delta(F)$ of $\Delta^{k+1}$.\\
    For every $j\in \set{2,...,n}\setminus S$, define
    \[
        E_j(F)=
        \begin{cases}
            \{0\},
            &\text{if }F\subseteq\mathbf L_{j,S},\\
            \{1\},
            &\text{if }F\subseteq\mathbf U_{j,S},\\
            I,
            &\text{if $F$ is contained in neither.}
        \end{cases}
    \]
    This is well-defined. Indeed, $F$ cannot be contained simultaneously in $\mathbf L_{j,S}$ and $\mathbf U_{j,S}$, since that would force
    \[
        x_j=\epsilon\lambda_S
        =
        1-\epsilon\lambda_S.
    \]
    for any $x \in F$, which is impossible since $\epsilon < \frac{1}{2}$.
    We have therefore defined a map from the set of nonempty faces of $\P$ to the set of nonempty faces of $\Delta^{k+1}\times I^m$ by
    \[
        F\longmapsto
        \left(\pi_\Delta(F),(E_j(F))_{j\in \set{2,...,n}\setminus S}\right).
    \]
    We now show that this map is a bijection.
    Let
    \[
        \left(D,(E_j)_{j\in \set{2,...,n}\setminus S}\right)
    \]
    be an arbitrary tuple of the above form. Consider the intersection of the $\mathbf{B}$-facets corresponding to the facets of $\Delta^{k+1}$ that contain $D$, together with $\mathbf L_{j,S}$ whenever $E_j=\{0\}$, and $\mathbf U_{j,S}$ whenever $E_j=\{1\}$. Since this is an intersection of facets, it is a face provided that it is nonempty.
    To see that it is nonempty, choose a point $q$ in the relative interior of the face of $\pi_\Delta(\P)$ corresponding to $D$. For every $j\in \set{2,...,n}\setminus S$, set
    \[
        x_j=
        \begin{cases}
            \epsilon\lambda_S(q),
            &\text{if }E_j=\{0\},\\
            1-\epsilon\lambda_S(q),
            &\text{if }E_j=\{1\},\\
            \frac12,
            &\text{if }E_j=I.
        \end{cases}
    \]
    The resulting point $(q,x_{j_1}, \dots, x_{j_m})$ belongs to exactly the facets prescribed by the tuple.
    It follows that every tuple determines a nonempty face of $\P$. Hence the map is surjective. It is injective as well: if two faces determine the same tuple, then they are contained in exactly the same $\mathbf{B}$-, $\mathbf{U}$-, and $\mathbf{L}$-facets. By \eqref{eqn:face_facets}, the two faces are equal.\\
    It remains to compare the partial orders. Suppose $F\subseteq F'$. Every facet containing $F'$ also contains $F$. Therefore,
    \[
        \pi_\Delta(F)\subseteq \pi_\Delta(F')
        \qquad\text{and}\qquad
        E_j(F)\subseteq E_j(F')
        \quad\text{for every }j\in \set{2,...,n}\setminus S.
    \]
    Conversely, if these coordinatewise containments hold, then every facet equality that holds for $F'$ also holds for $F$. Therefore $F\subseteq F'$.
    We have thus obtained an isomorphism of their complete face lattices (after extending to the empty face in the obvious way).
\end{proof}

\section{A regular prismatic subdivision}
\label{sec:regular-subdivision}

We now show that the prisms $\P$ are all the maximal cells of a regular prismatic subdivision of $\G$ (Definition \ref{def:polytopal-subdivision}, Theorem \ref{thm:regular-subdivision}).
\begin{definition}
    Define the height function
    \[
    h(v) = \begin{cases}
        0, \quad \text{if }v \text{ is a cube point}\\
        1, \quad \text{if }v \text{ is an added point}
    \end{cases}
    \]
    and the \textit{lifted Goodwillie polytope}
    \[
    \left({G}_n^\epsilon\right)^h = \conv \left\{ (v,h(v)): v\text{ is a defining point of }\G \right\} \subset \R^{n+1}.
    \]
\end{definition}
We have that $\left( \G \right)^h$ is an $(n+1)$-dimensional convex polytope. 

Denote the last coordinate of $\R^{n+1}$ by $z$. We will show that the planes 
\[
z-\ell_S(x)=0 , \quad S\subseteq \set{2,\ldots,n}
\]
where
\[
\quad \ell_S(x) = - \frac{\sum_{i\in S} x_i}{1-(\abs{S}-1)\epsilon}
\]
are the supporting planes of the lower facets of $\left( \G \right)^h$, which are in bijection with the maximal cells of the regular subdivision of $\G$ induced by $h$ (Theorem \ref{thm:regular-subdivision}).

\begin{proposition}\label{prop:prisms-regular-cells}
    For every $S\subseteq \set{2,...,n}$, the hyperplane $H_S:z = \ell_S(x)$ supports $\left( \G \right)^h$. The lower face $H_S \cap \left( \G \right)^h$ projects bijectively onto $\P$ under $\rho$. Thus, every $\P$ is a maximal cell of the regular subdivision of $\G$ induced by $h$.
\end{proposition}
\begin{proof}
    Let $k=\abs{S}$. Let $v$ be a vertex of $\left( \G \right)^h$, then $v$ is of the form
    \[
      (c_{t,A}, h(c_{t,A})) = (c_{t,A}, 0)       \quad \text{or} \quad  (u_{r,B}^\epsilon,h(u_{r,B}^\epsilon))=(u_{r,B}^\epsilon,1).  
    \]
    In the first case, $0 - \ell_S(c_{t,A}) = \frac{\sum_{i\in S} (c_{t,A})_i}{1-(k-1)\epsilon}\geq 0$, and equality holds precisely when $A\cap S = \emptyset$, i.e.~when $c_{t,A}$ is a vertex of $\P$. In the second case, $1 - \ell_S(u_{r,B}^\epsilon) = 1+ \frac{\sum_{i\in S} (u_{r,B}^\epsilon)_i}{1-(k-1)\epsilon}$. Here we must consider two subcases; if $r \notin S$ then 
    \[ 1+ \frac{\sum_{i\in S} (u_{r,B}^\epsilon)_i}{1-(k-1)\epsilon} \geq 0  \]
    is clear, as \(\sum_{i\in S} (u_{r,B}^\epsilon)_i\) is a sum of positive real numbers. If $r \in S$ , then
    \begin{align*}   
    1-(k-1)\epsilon+ \sum_{i\in S} (u_{r,B}^\epsilon)_i &= -(k-1)\epsilon + \sum_{i\in S\setminus\{r\}} (u_{r,B}^\epsilon)_i \\
    &= -(k-1)\epsilon + \sum\limits_{i \in S\cap B} (1-\epsilon) +\sum\limits_{i\in S\setminus (B\cup\{r\}) }\epsilon \\
    &\geq -(k-1)\epsilon  +\sum\limits_{i\in S\setminus \{r\} }\epsilon = -(k-1)\epsilon  + (k-1)\epsilon = 0,
    \end{align*}
    with equality precisely when $r \in S$ and $B\cap S = \emptyset$, i.e.~when $u_{r,B}^\epsilon$ is a vertex of $\P$.\\ 
    Therefore $\conv(\{ (v,h(v)) \, \vert \, v \in V(\P)   \})$ is a lower face of $(\G)^h$ (\ref{def:lower-face}) and, since $\dim((\G)^h) = n+1$, by Theorem \ref{thm:regular-subdivision}, $\P$ is a maximal cell of the polytopal subdivision of $\G$ induced by $h$.
\end{proof}

\begin{corollary}\label{cor:prism-intersection}
    For $S, T\subseteq \set{2,...,n}$, the intersection of $\P$ and $P_T^\epsilon$ is a common face, satisfying
    \begin{enumerate}
        \item $\P \cap P_T^\epsilon = \conv \left(  \set{ c_{t,A} \ \vert \ A\cap (S\cup T)=\emptyset } \cup \set{u^\epsilon_{r,B} \ \vert \ r\in S\cap T, B\cap(S\cup T)= \emptyset} \right),$
        \item $\P \cap P_T^\epsilon \cong \Delta^{\abs{S\cap T}+1} \times I^{n-1-\abs{S\cup T}}\ ,$ and
        \item $\dim (\P \cap P_T^\epsilon) = n- \abs{(S \cup T) \setminus (S \cap T)}$.
    \end{enumerate}
    In particular, whenever $\P , P_T^\epsilon$ share a facet we have $\abs{(S \cup T) \setminus (S \cap T)} = 1$, so that $T=S\cup \set{j}, j\notin S$ (or vice-versa), and 
    \[
    \P \cap P_{S\cup \set{j}}^\epsilon  = \mathbf{L}_{j,S} = \mathbf{B}_{\upsilon_j, S\cup \set{j}}.
    \]
\end{corollary}
\begin{proof}
By the preceding proposition and
Theorem \ref{thm:regular-subdivision},
$P_S^\epsilon$ and $P_T^\epsilon$ are cells of the same
regular subdivision, so their intersection is a common face.
It is nonempty since both contain $c_{0,\emptyset}$.

The vertices of a common face are precisely the common vertices
of the two polytopes. By Proposition \ref{prop:prism_facets}
and Definition \ref{def:P}, these are
\[
c_{t,A}\quad\text{with }A\cap(S\cup T)=\emptyset,
\]
and
\[
u^\epsilon_{r,B}\quad\text{with }
r\in S\cap T,\qquad B\cap(S\cup T)=\emptyset.
\]
Taking their convex hull proves (1).

We now describe this face in the product structure of
$P_S^\epsilon$ given by Theorem \ref{thm:prism-comb-type},
\[\P \cong  \Delta^{\abs{S}+1}\times I^{n-\abs{S}-1}, \]
which partitions the coordinates of $\R^n$ into simplex and cube coordinates.
For each $i\in S\setminus T$, impose the facet
$\mathbf B_{\upsilon_i,S}$, deleting the simplex vertex
$\upsilon_i$. For each $j\in T\setminus S$, impose
$\mathbf L_{j,S}$, fixing the corresponding cube coordinate
at its lower value $\epsilon \lambda_{S}(q)$. The vertices satisfying the above are exactly those listed in (1), so these
conditions describe $P_S^\epsilon\cap P_T^\epsilon$.

The remaining simplex vertices are
\[
\zeta_0,\ \zeta_1,\quad
\upsilon_r\quad(r\in S\cap T),
\]
and the free cube coordinates are indexed by
$\set{2,...,n}\setminus(S\cup T)$. Consequently,
\[
P_S^\epsilon\cap P_T^\epsilon
\cong
\Delta^{|S\cap T|+1}\times I^{n-1-|S\cup T|},
\]
which proves (2). Its dimension is
\[
(|S\cap T|+1)+(n-1-|S\cup T|)
=n-|(S \cup T) \setminus (S \cap T)|,
\]
proving (3).

Since both prisms have dimension $n$, their intersection
is a facet of both precisely when $|(S \cup T) \setminus (S \cap T)|=1$.
After interchanging $S$ and $T$ if necessary, write
$T=S\cup\{j\}$. The face description above imposes only
$\mathbf L_{j,S}$ in $P_S^\epsilon$, and only
$\mathbf B_{\upsilon_j,S\cup\{j\}}$ in
$P_{S\cup\{j\}}^\epsilon$. Hence
\[
P_S^\epsilon\cap P_{S\cup\{j\}}^\epsilon
=
\mathbf L_{j,S}
=
\mathbf B_{\upsilon_j,S\cup\{j\}}.
\]
\end{proof}
We will make use of the following technical proposition, which allows us to pullback the equations for the defining planes of $\P$ in $\fancy{H} \times \R^{n-k-1}$ to ``ambient coordinates'' in $\R^n$.
\begin{lemma}\label{lem:prism-to-ambient}
Let $S\subseteq \set{2,...,n}$ and $x\in P_S^\epsilon$. Write
$q=\pi_\Delta(x)$, and let $z_0,z_1,(w_i)_{i\in S}$
be the barycentric coordinates of $q$. Then
\begin{align}
&\lambda_S(q)
=-\frac{\sum_{i\in S}x_i}
         {1-(\abs{S}-1)\epsilon},\label{eq:prism-to-ambient1}\\
&w_i
=\frac{\epsilon\lambda_S(q)-x_i}{1+\epsilon},
\qquad i\in S,\label{eq:prism-to-ambient2}\\
&z_1
=x_1-\sum_{i\in S}\frac{i-1}{n}w_i,\label{eq:prism-to-ambient3}\\
&z_0 =1-z_1-\lambda_S(q).\label{eq:prism-to-ambient4}
\end{align}
\end{lemma}

\begin{proof}
Since $x\in P_S^\epsilon$, we can
write $x$ as a convex combination of the defining points
of $P_S^\epsilon$. Applying $\pi_\Delta$ to this combination
gives the barycentric expansion
\[
q=z_0\zeta_0+z_1\zeta_1+\sum_{r\in S}w_r\upsilon_r.
\]
By uniqueness of barycentric coordinates, $z_0$ is the
total weight assigned to cube points with first coordinate
$0$, $z_1$ is the total weight assigned to cube points with first coordinate $1$, and $w_r$ is the total weight assigned to added
points of type $r$.

Fix $i\in S$. Taking the $i$-th coordinate on both sides
of the convex combination defining $x$ expresses $x_i$
as the same weighted sum of the $i$-th coordinates of
the defining points. These coordinates are zero for
cube points, $-1$ for added points of type $i$, and
$\epsilon$ for added points of every other type.
Indeed, the defining added points satisfy
$B\cap S=\emptyset$, so their coordinates indexed by
$S$ cannot take the high value $1-\epsilon$.
Grouping the contributions by type therefore gives
\[
x_i
=0(z_0+z_1)-w_i
+\epsilon\sum_{r\in S\setminus\{i\}}w_r.
\]
Since $\lambda_S(q)=\sum_{r\in S}w_r$, this becomes
\[
x_i
=-w_i+\epsilon\bigl(\lambda_S(q)-w_i\bigr)
=\epsilon\lambda_S(q)-(1+\epsilon)w_i.
\]
Solving for $w_i$ proves \eqref{eq:prism-to-ambient2}.

Summing these equalities over $i\in S$ gives
\[
\begin{aligned}
\sum_{i\in S}x_i
&=\abs{S}\epsilon\lambda_S(q)
  -(1+\epsilon)\sum_{i\in S}w_i\\
&=-\bigl(1-(\abs{S}-1)\epsilon\bigr)\lambda_S(q).
\end{aligned}
\]
The factor $1-(\abs{S}-1)\epsilon$ is positive in the
assumed parameter range, so division proves the first
identity \eqref{eq:prism-to-ambient1}.

Similarly, taking the first coordinate of the convex
combination defining $x$, cube points contribute
$0z_0+1z_1$, while added points of type $i$ contribute
$\frac{i-1}{n}w_i$. Hence
\[
x_1=z_1+\sum_{i\in S}\frac{i-1}{n}w_i,
\]
which proves the third identity. Finally, the barycentric
coordinates sum to one:
\[
z_0+z_1+\sum_{i\in S}w_i=1.
\]
Using $\sum_{i\in S}w_i=\lambda_S(q)$ gives the fourth
identity \eqref{eq:prism-to-ambient4}.
\end{proof}

\begin{proposition}\label{prop:external-facets}
For every $S\subseteq \set{2,...,n}$, the facets
$\mathbf B_{\zeta_0,S}$, $\mathbf B_{\zeta_1,S}$, and
$\mathbf U_{j,S}$, for $j\in \set{2,...,n}\setminus S$, of $P_S^\epsilon$
are also facets of $G_n^\epsilon$. Write
\[
  \mathbf C_{t,S}:=\mathbf B_{\zeta_{1-t},S},
  \qquad t\in\{0,1\}.
\]
Among the defining points of $G_n^\epsilon$:
\begin{enumerate}
\item $\mathbf C_{t,S}$ contains exactly the cube points
$c_{t,D}$ and the added points $u_{r,D}^\epsilon$ with
$r\in S$, both subject to $D\cap S=\emptyset$.
\item $\mathbf U_{j,S}$ contains exactly the cube points
$c_{t,D}$ of either type $t\in\{0,1\}$ and the added points
$u_{r,D}^\epsilon$ with $r\in S$, both subject to
$j\in D$ and $D\cap S=\emptyset$.
\end{enumerate}
\end{proposition}
\begin{proof}
The local facets in the statement are defined by the equalities
$z_1=0$, $z_0=0$, and $1-x_j-\epsilon\lambda_S=0$.
We will show that the corresponding inequalities hold on the whole polytope,
with equality precisely at the defining points listed above.

The formulas in Lemma \ref{lem:prism-to-ambient} extend
$\lambda_S,z_0,z_1$ to affine functions of $x\in\R^n$;
we write $\lambda_S(x),z_0(x),z_1(x)$ for these extensions.
In particular,
\[
\lambda_S(x)=-\frac{\sum_{i\in S}x_i}{1-(\abs{S}-1)\epsilon}.
\]
Substituting the expressions for $w_i$ into those for $z_1$ and $z_0$ gives
\begin{align*}
z_1(x)&=x_1+\frac{1}{1+\epsilon}\left[
\sum_{i\in S}\frac{i-1}{n}x_i
+\frac{\epsilon\sum_{i\in S}\frac{i-1}{n}}
       {1-(\abs{S}-1)\epsilon}\sum_{i\in S}x_i\right],\\
z_0(x)&=1-x_1+\frac{1}{1+\epsilon}\left[
\sum_{i\in S}\left(1-\frac{i-1}{n}\right)x_i
+\frac{\epsilon\sum_{i\in S}\left(1-\frac{i-1}{n}\right)}
       {1-(\abs{S}-1)\epsilon}\sum_{i\in S}x_i\right].
\end{align*}
The denominator $1-(\abs{S}-1)\epsilon$ is positive.
Since $0<(i-1)/n<1$ for $i\in S$, the coefficient of
every $x_i$, $i\in S$, in each expression is strictly positive.

For a cube point $c_{t,D}$, the coordinates indexed by $S$ are
nonnegative. Thus $z_1(c_{t,D})\geq0$, with equality exactly
when $t=0$ and $D\cap S=\emptyset$. Likewise,
$z_0(c_{t,D})\geq0$, with equality exactly when $t=1$
and $D\cap S=\emptyset$.

Consider an added point $u_{r,D}^\epsilon$. If $r\notin S$,
its coordinates indexed by $S$ are positive and its first
coordinate lies strictly between $0$ and $1$. Hence both
$z_0$ and $z_1$ are strictly positive there.
If $r\in S$, consider first $u_{r,D\setminus S}^\epsilon$.
This is a defining point of $P_S^\epsilon$ whose projection is
$\upsilon_r$, so $z_0=z_1=0$ there. We can obtain
$u_{r,D}^\epsilon$ from $u_{r,D\setminus S}^\epsilon$ by increasing each coordinate indexed by
$D\cap S$ from $\epsilon$ to $1-\epsilon$, leaving all
other coordinates fixed. Each increase is positive and has a
strictly positive coefficient in both affine functions.
Consequently, $z_0,z_1\geq0$ at $u_{r,D}^\epsilon$,
with equality exactly when $D\cap S=\emptyset$.

For $j\notin S$, the remaining supporting function is
\[
H_{U_j,S}(x):=1-x_j-\epsilon\lambda_S(x)
=1-x_j+\frac{\epsilon\sum_{i\in S}x_i}
                  {1-(\abs{S}-1)\epsilon}.
\]
At a cube point, its value is
\[
1-(c_{t,D})_j+
\frac{\epsilon\abs{D\cap S}}{1-(\abs{S}-1)\epsilon},
\]
which is nonnegative, with equality exactly when
$j\in D$ and $D\cap S=\emptyset$.
At an added point with $r\notin S$, it is strictly positive:
$(u_{r,D}^\epsilon)_j<1$ and the sum over $S$ is nonnegative.
Finally, if $r\in S$, then
\[
\sum_{i\in S}(u_{r,D}^\epsilon)_i
=-1+(\abs{S}-1)\epsilon+(1-2\epsilon)\abs{D\cap S},
\]
so
\[
H_{U_j,S}(u_{r,D}^\epsilon)
=1-(u_{r,D}^\epsilon)_j-\epsilon
+\frac{\epsilon(1-2\epsilon)\abs{D\cap S}}
       {1-(\abs{S}-1)\epsilon}.
\]
Here $j\neq r$, so the first three terms sum to $0$ if
$j\in D$ and to $1-2\epsilon>0$ otherwise. This proves
nonnegativity and the remaining stated incidence conditions.

Since these functions are affine, their nonnegativity on the
defining points implies nonnegativity on $G_n^\epsilon$.
Each zero set cuts out the convex hull of the defining points
where equality holds. By the stated incidence conditions, this is exactly the
corresponding local facet of $P_S^\epsilon$.
It has dimension $n-1$, so it is also a facet of $G_n^\epsilon$.
\end{proof}
\begin{remark}
    Note that in the above, $\mathbf{C}_{0,S}$ corresponds to $\mathbf{B}_{\zeta_1,S}$. The reason is that the $\mathbf{B}_{\zeta_1,S}$ corresponds to \textbf{the face opposite to $\zeta_1$}. This unfortunate notational detail is important to state Theorem \ref{thm:self-duality} cleanly. 
\end{remark}
Let's now summarize what we have learned about the facets of $\P$.
\begin{proposition}\label{prop:facets-of-P}
    For $S\subseteq \set{2, \dots, n}$, all facets of $\P$ fall into one of the following cases:
    \begin{enumerate}
        \item For each $i\in S$, the facet $\mathbf{B}_{\upsilon_i,S}$ of $\P$ is equal to the facet $\mathbf{L}_{i,S\setminus \set{i}}$ of $P_{S\setminus \set{i}}^\epsilon$.
        \item For each $j\in \set{2,...,n}\setminus S$, the facet $\mathbf{L}_{j,S}$ of $\P$ is equal to the facet $\mathbf{B}_{\upsilon_j, S\cup \set{j}}$ of $P_{S\cup \set{j}}^\epsilon$.
        \item The two facets $\mathbf{B}_{\zeta_0, S}, \mathbf{B}_{\zeta_1,S}$, and all $\mathbf{U}$-facets are facets of $\G$.
    \end{enumerate}
\end{proposition}
\begin{proof}
The complete list of facets is given in
Corollary \ref{cor:facets_P}. The identifications in
(1) and (2) follow from Corollary \ref{cor:prism-intersection},
and (3) is Proposition \ref{prop:external-facets}.
\end{proof}
We can now state the main theorem of the section. 
\begin{theorem}\label{thm:subdivision}
    Let $n\geq 3, 0<\epsilon< \min \{ \frac{1}{2} , \frac{1}{n-2} \}$. For every $S \subseteq \set{ 2, \dots, n }$, the polytope $\P$ is $n$-dimensional and
    \[
    \P \iso \Delta^{\abs{S}+1} \times I^{n-\abs{S}-1}. 
    \]
    The $\P$ are precisely the maximal cells of a regular prismatic subdivision of $\G$.\\ 
    In particular, \[ \G = \bigcup_{S\subseteq \set{2, \dots, n}} \P ,\]
    and two cells $\P$ and $P_T^\epsilon$ share a facet if and only if their labels differ by exactly one element.\\
    Moreover, the vertices of $G_n^\epsilon$ are precisely
its defining points $c_{t,A}$ and $u^\epsilon_{r,B}$.
\end{theorem}
\begin{proof}
By Proposition \ref{prop:prisms-regular-cells}, every
$P_S^\epsilon$ is a maximal cell of the regular subdivision
induced by $h$. Its combinatorial type is given by
Theorem \ref{thm:prism-comb-type}, and the assertion about
shared facets follows from
Corollary \ref{cor:prism-intersection}.
    It remains to show that there are no additional maximal cells.
The maximal cells of a subdivision of a convex polytope are
connected through shared facets.
If an additional maximal cell existed, there would be a maximal cell
not among the listed prisms $P_S^\epsilon$ sharing a facet with one of them.
By Proposition \ref{prop:facets-of-P}, however, every facet
of $P_S^\epsilon$ either lies on the boundary of $\G$
or is already shared with another listed prism, so no additional maximal cell can share a facet with a prism.

Thus the prisms are all the maximal cells, and their union
is $G_n^\epsilon$ by Theorem \ref{thm:regular-subdivision}.\\
Finally, every cube point $c_{t,A}$ is a vertex of the
facet $\mathbf C_{t,\emptyset}$, and every added point
$u^\epsilon_{r,B}$ is a vertex of the facet
$\mathbf C_{0,\{r\}}$, by
Propositions \ref{prop:prism_facets} and
\ref{prop:external-facets}.
A vertex of a face is a vertex of the whole polytope,
so all defining points are vertices of $G_n^\epsilon$.
Conversely, every vertex of the convex hull of a finite
set belongs to that set.
\end{proof}
\begin{figure}
    \centering \includegraphics[width=\textwidth]{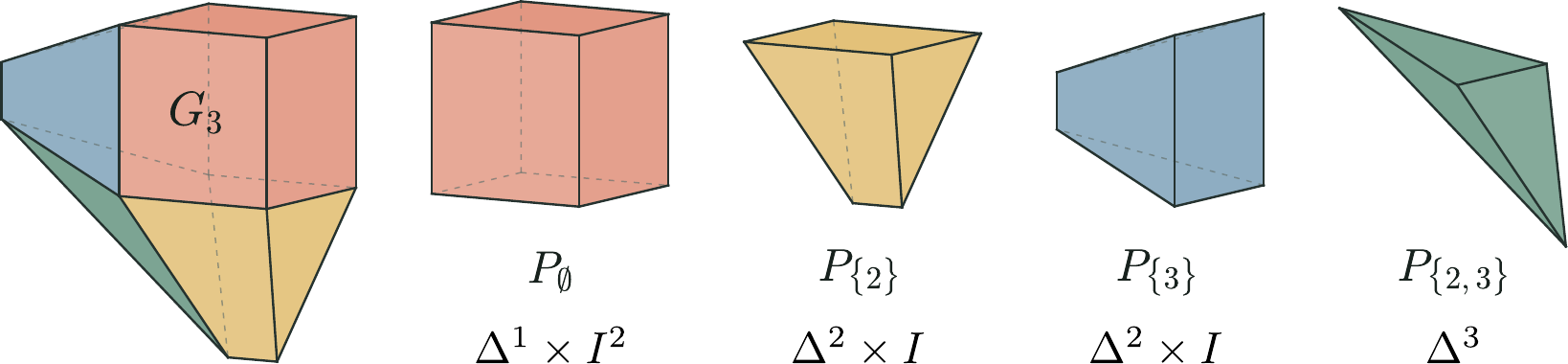}
    \caption{The regular prismatic subdivision of $G_3^\epsilon$ from Theorem~\ref{thm:subdivision}. Its four maximal cells, shown separately on the right, are combinatorially equivalent to $\Delta^1\times I^2$, two copies of $\Delta^2\times I$, and $\Delta^3$. Matching colors identify the cells in the assembled polytope. The superscript $\epsilon$ is omitted from the figure labels.}
    \label{fig:g3-decomp}
\end{figure}
\begin{corollary}\label{cor:stable-range}
    Let $n\geq 3$. If $\epsilon,\epsilon' \in \left( 0,\min\set{\frac{1}{2},\frac{1}{n-2}}\right) $ then $\G$ and $G_{n}^{\epsilon'}$ are combinatorially equivalent.
\end{corollary}
From now on, we fix $\epsilon$ in the admissible range and suppress
it from the notation, writing $G_n$, $P_S$, and $u_{r,B}$ for
$G_n^\epsilon$ , $P_S^\epsilon$, and $u_{r,B}^\epsilon$. We use $G_n$ for both this
realization and its combinatorial type, as determined by context.

An immediate consequence of the prismatic decomposition of $G_n$ is statement (2) in Theorem \ref{theorem1}, which we record in the next corollary.

\begin{corollary}\label{Gisprismatic} Each $G_n$ is a prismatic $n$-dimensional polytope.
\end{corollary}

\section{Combinatorial self-duality}\label{sec:self-duality}
A polytope $P$ is said to be \textit{combinatorially self-dual} if there is an isomorphism of posets between its face lattice $\mathcal F(P)$ and its opposite $\mathcal F(P)^{op}$; equivalently, if there is a bijection $\Phi: \mathcal F(P)\rightarrow \mathcal F(P)$ such that, for all $F,F'\in \mathcal F(P)$,  $F\subseteq F'$ if and only if $ \Phi(F') \subseteq \Phi(F)$. 
In this section, we prove statement (1) in Theorem \ref{theorem1}: each $G_n$ is combinatorially self-dual. 

Recall from Theorem \ref{thm:subdivision} that $G_n$ has vertex set
    \[
    V(G_n) = \set{c_{t,A} \ \vert \ t\in \set{0,1}, A \subseteq \set{2,...,n}} \sqcup \set{u_{r,B} \ \vert \ r\in \set{2,...,n}, B\subseteq \set{2,...,n}\setminus \set{r}},
    \]
    and, by Theorem \ref{thm:subdivision} and Proposition \ref{prop:facets-of-P}, its facet set is
    \[
    \begin{aligned}
    \text{Facets}(G_n)
    &= \set{\mathbf C_{t,A} \ \vert \ t\in \set{0,1}, A \subseteq \set{2,...,n}}\\
    &\quad{}\sqcup \set{\mathbf U_{r,B} \ \vert \ r\in \set{2,...,n}, B\subseteq \set{2,...,n}\setminus \set{r}}.
    \end{aligned}
    \]

The map of sets 
\begin{align*}
    \Phi: V(G_n) &\longrightarrow \text{Facets}(G_n)\\
    c_{t,A} &\longmapsto \mathbf C_{t,A} \\
    u_{r,B} 
    &\longmapsto \mathbf U_{r,B}
\end{align*}
is clearly a bijection.
We can extend this map to the entire face lattice of $G_n$ by
\begin{align}
    \Phi: \mathcal F(G_n) &\longrightarrow \mathcal F(G_n) \\
    F &\longmapsto \bigcap\limits_{v \in V(F)} \Phi(v)
\end{align}
We set $\Phi(\emptyset)=G_n$. Since an intersection of facets of a polytope determines a (possibly empty) face, $\Phi$ is well defined. Order reversal follows immediately from definition.
\begin{proposition}
    If $F \subseteq F'$ then $\Phi(F') \subseteq \Phi(F)$.
\end{proposition}
We now state two technical propositions.
\begin{proposition}\label{lem:vertex-facet-incidence}
    For all vertices $v,w$ of $G_n$
    $ v \in \Phi(w)$ if and only if $w \in \Phi(v)$.
\end{proposition}
\begin{proof}
The equality cases in Proposition \ref{prop:external-facets} give the following three cases.
For two cube vertices, we have that 
$c_{t,A}\in\Phi(c_{s,D})$
if and only if
$t=s$ and $A\cap D=\emptyset,$
which is symmetric in the two vertices.\\
For a cube vertex and an added vertex, we have that 
$c_{t,A}\in\Phi(u_{r,B})$  if and only if 
$r\in A\text{ and }A\cap B=\emptyset$
if and only if 
$u_{r,B}\in\Phi(c_{t,A})$.

Finally, for two added vertices,
$u_{r,B}\in\Phi(u_{s,D})$ if and only if
\[r\in D,\quad s\in B,\text{ and }B\cap D=\emptyset.
\]
This condition is again symmetric in the two vertices.
\end{proof}
For a face $F$, denote \[
\rchi(F)= \set{ H \in \text{Facets}(G_n) \ \vert \ F \subseteq H }.
\]
\begin{proposition}\label{lem:dual-face-vertices}
    For every face $F$ of $G_n$,
    \[
    V(\Phi(F)) = \Phi^{-1}(\rchi (F)) \ .
    \]
\end{proposition}
\begin{proof}
    Let $w \in V(G_n)$ and $F \in \mathcal F(G_n)$. We have the following equivalences.
    \begin{align*}
        w\in \Phi(F) &\iff w \in \Phi(v) \quad \text{ for all } v \in V(F) \quad \text{by definition}\\
        &\iff v \in \Phi(w) \quad \text{ for all } v \in V(F) \quad  \text{by } \ref{lem:vertex-facet-incidence}\\
        &\iff F \subseteq \Phi(w)\\
        &\iff \Phi(w) \in \rchi(F)\\
        &\iff w \in \Phi^{-1}(\rchi(F))  \ .
    \end{align*}
\end{proof}
\begin{theorem}\label{thm:self-duality}
    The map
    \[
    \Phi : \mathcal F (G_n) \rightarrow \mathcal F(G_n)
    \]
    is an order-reversing involution. In particular, $G_n$ is combinatorially self-dual.
\end{theorem}
\begin{proof}
    We have already shown order-reversal, so we need only show that $\Phi$ is an involution, i.e.~ $\Phi^2 = \id$, which follows from the previous proposition:
    \begin{equation*}
        \Phi^2(F) = \Phi\left(\bigcap_{v\in V(F)} \Phi(v)\right)
        \ =\bigcap_{w \in V(\Phi(F))} \Phi(w) \ =\bigcap_{w \in \Phi^{-1}(\rchi(F))} \Phi(w) =\bigcap_{H \in \rchi(F)} H = F \ .
    \end{equation*}
\end{proof}
\begin{corollary}\label{cor:palindromic-f-vector}
    The $f$-vector of $G_n$ is a palindrome, i.e.
    \[
    f_i(G_n) = f_{n-1-i}(G_n) \quad \text{for } -1\leq i\leq n .
    \]
\end{corollary}

Recall that one may produce new self-dual polytopes by taking joins and elongated pyramids. If $P$ and $Q$ are two convex polytopes their \textit{join}, denoted by $P \star Q$, is defined by positioning $P$ and $Q$ into skew affine subspaces and taking the convex hull. The resulting polytope has dimension $\text{dim}(P)+ \text{dim}(Q) + 1$. If $Q$ is a point $q$, $P \star q$ is also called the \textit{pyramid} over $P$.  If $P$ and $Q$ are self-dual, then so is $P \star Q$.

The \textit{elongated pyramid} over $P$ is the $(\text{dim}(P)+1)$-dimensional polytope defined by taking the union of $P \times I$ and $(P \times \{1\}) \star q$, where $q$ is a $0$-dimensional polytope. If $P$ is self-dual so is the elongated pyramid over $P$. 

We show that the Goodwillie polytopes do not arise through these constructions for $n \geq 3$. This is statement (3) in Theorem \ref{theorem1}.

\begin{proposition}
     For $n\geq 3$, $G_n$ is not combinatorially equivalent to a join of two non-empty polytopes or to an elongated pyramid over an $(n-1)$-dimensional polytope. 
\end{proposition}
\begin{proof}
We first show that there are no (non-empty) polytopes $P$ and $Q$ such that $G_n \cong P \star Q$. By Corollary \ref{Gisprismatic}, every facet of $G_n$ is a prism. Since every prism $\Delta^k \times I^l$ is a simple polytope (an $n$-polytope in which every vertex is adjacent to exactly $n$-edges) being a product of two simple polytopes, every facet of $G_n$ is a simple polytope. 

We claim that if $D$ is a simple $m$-polytope with $d=|V(D)|$ vertices such that $D=A \star B$ for two polytopes $A$ and $B$, then $D$ must be a simplex. Suppose $A$ and $B$ are of dimension $p$ and $q$, respectively, so $p+q=m-1$. A vertex $v\in A$ in $D$ is adjacent to every vertex of $B$, and to at least $p$ vertices of $A$. So $\text{deg } v\geq |V(B)|+p$. Similarly, $\text{deg } w\geq |V(A)|+q$ for any vertex $w\in B$ in $D$. By simplicity $\text{deg } v= m = \text{deg } w$, so adding the two inequalities gives $2m\geq d+m-1$. Hence $m\ge d-1$. Since $d\geq m+1$, we get $d=m+1$, which implies $D \cong \Delta^m$. 

Suppose $G_n \cong P \star Q$ and the dimension of both $P$ and $Q$ is at least $1$. Let $F$ be a facet of $P$.  Then $F \star Q$ is a facet of $G_n$, so by the above claim it must be a simplex and consequently $Q$ must be a simplex as well. Similarly, $P$ must be a simplex. This means that $G_n$ must be a simplex, which is not the case for $n \geq 2$.

Suppose one of $P$ and $Q$, say $P$, is a point so that $G_n$ is a pyramid over $Q$ with cone point $P$. For any facet $F$ of $Q$, we have facet $P \star F$ of $G_n$. By the above claim, $P \star F$ must be a simplex, which implies $F$ must be a simplex as well, so all facets of $Q$ must be simplices. Since $Q$ is a facet of $G_n$, which must be a prism of dimension $n-1$, we must have either $Q \cong \Delta^{n-1}$ or $Q \cong I^2$. If $Q \cong \Delta^{n-1}$ then $G_n \cong \Delta^n$, which is not true for $n\geq 2$.  If $Q \cong I^2$, then $n=3$, but in this case we can check explicitly that $G_3$ is not a pyramid over a square. 

The polygon $G_2$ is a pentagon, which is an elongated pyramid over $\Delta^1$. For $n \geq 3$, $G_n $ has an even number of vertices while an elongated pyramid over a positive-dimensional polytope $P$ has $2|V(P)|+1$ vertices.

\end{proof}

\section{Face classification and the \texorpdfstring{$f$-polynomial}{f-polynomial}}
\label{sec:face-enumeration}
For $0\leq d\leq n-1$, let $f_d(G_n) $ be the number of
$d$-dimensional faces of $G_n$ and define the $f$-polynomial of $G_n$ by
\begin{equation}
\label{eq:face-polynomial-definition}
  \mathfrak f_n(x)
  \coloneqq
  \sum_{d=0}^{n-1}f_d(G_n)x^d.
\end{equation}
Thus a $d$-dimensional face contributes $x^d$.  The empty face and
$G_n$ itself are not included; as usual, we set
$f_{-1}=f_n=1$ when discussing the complete $f$-vector.

In this section, we obtain a closed formula for the $f$-polynomial
by classifying the proper nonempty faces of $G_n$ and counting
them by dimension. To count each face exactly once, we need
a unique description of it.

Our starting point is the vertex--facet incidence rules.
A vertex belongs to a face precisely when it belongs to every
facet containing that face. Applying these rules to all the
facets containing a face $F$ determines which cube and added
vertices belong to $F$. We describe this vertex set by a label $(U,A,L,T)$,
where $U,A,L$ are pairwise disjoint subsets of $\{2,\ldots,n\}$ and
$T\subseteq\{0,1\}$. These sets record the allowed vertex types and the
restrictions on their coordinates. We then characterize the
possible labels and show that they are in bijection with the
proper nonempty faces, and the dimension and combinatorial
type of each face is determined by its label
(Theorem~\ref{thm:canonical-face-classification}).

\subsection{Canonical face labels}

For a face $F$, let
\[
  \rchi(F)
  \coloneq
  \{H\in\Facets(G_n):F\subseteq H\}.
\]
If $F$ is proper, then $\rchi(F)\neq\emptyset$, and $F$ is the
intersection of its containing facets
\begin{equation}
\label{eq:face-as-all-containing-facets}
  F=\bigcap_{H\in\rchi(F)}H.
\end{equation}

To describe the vertices of $F$, we apply the incidence rules
of Proposition~\ref{prop:external-facets} to every facet
containing $F$.
Recall what the subset $D$ in the subscript of a vertex records:  for a cube vertex
$c_{t,D}$, $D$ is the set of indices $i\in \set{2,...,n}$ for which the $i$-th
coordinate is $1$, rather than $0$.  For an added vertex
$u_{a,D}$, $D$ is the set of indices $i\in \set{2,...,n}\setminus\{a\}$ for
which the $i$-th coordinate is $1-\epsilon$, rather than $\epsilon$.  In either
case, we call $D$ the \emph{high-coordinate set}.

We use the first entry in the subscript of a vertex to record its \emph{type}:
the cube vertex $c_{t,D}$ has type $t\in\{0,1\}$, and the added vertex
$u_{a,D}$ has type $a\in \set{2,...,n}$. The type determines the image of a vertex under the projection to the simplex factor of any prism containing it (see \ref{eqn:proj} for the definition of these maps).

Recall also what the facet subscripts record.  In both $\mathbf C_{t,S}$ and
$\mathbf U_{j,S}$, the subset $S\subseteq \set{2,...,n}$ indexes the local prism
$P_S$ from which the facet was obtained; we call $S$ the \emph{support} of
the facet.  The index $t\in\{0,1\}$ in $\mathbf C_{t,S}$ specifies
the allowed cube-vertex type: every cube vertex belonging
to this facet has the form $c_{t,D}$.
The index $j$ in $\mathbf U_{j,S}$ specifies the coordinate
that must be high for a vertex to belong to the facet.

Thus the support $S$ of either facet excludes every index in $S$ from
the high-coordinate set.  For a cube vertex, those coordinates are therefore
low.  For an added vertex $u_{a,D}$, membership in either kind of facet
requires $a\in S$; coordinate $a$ has value $-1$, while
the other coordinates in $S$ are low.  A $\mathbf U_{j,S}$-facet additionally
forces coordinate $j$ to be high, while a $\mathbf C_{t,S}$-facet restricts the
cube type to $t$.

Let $F$ now be a proper nonempty face.  We read four pieces of information
from $\rchi(F)$:
\begin{itemize}
  \item $U$ is the set of indices $j$ for which a facet
  $\mathbf U_{j,S}$ belongs to $\rchi(F)$ for some $S$.
  These coordinates are forced high.

  \item $A$ is the intersection of the supports of all facets
  in $\rchi(F)$. These are the possible types of added vertices
  belonging to $F$.

  \item $L$ consists of the indices that belong neither to $U$
  nor to the support of any facet in $\rchi(F)$.  These coordinates remain
  free.

  \item $T\subseteq\{0,1\}$ records the cube types allowed by the
  containing $\mathbf C$-facets.  Thus $T=\{0,1\}$ if there is no containing
  $\mathbf C$-facet, $T=\{t\}$ if all of them have type $t$, and
  $T=\emptyset$ if both types occur.
\end{itemize}

\begin{lemma}\label{lem:face-vertices}
Let $F$ be a proper nonempty face of $G_n$, and let
$U,A,L,T$ be the sets defined above. Then $U,A,L$ are
pairwise disjoint, and
\begin{equation}
\label{eq:canonical-face-vertices}
\begin{split}
  V(F)
  ={}&
  \{c_{t,U\cup E}:t\in T,\ E\subseteq L\}\\
  &\mathbin{\sqcup}
  \{u_{a,U\cup E}:a\in A,\ E\subseteq L\}.
\end{split}
\end{equation}
\end{lemma}
\begin{proof}
Choose a vertex of $F$. By
Proposition~\ref{prop:external-facets}, its high-coordinate set
contains $U$ and avoids the support of every facet containing
$F$. Hence $U$ is disjoint from every such support.
Since $A$ is their intersection and $L$ is disjoint from both
$U$ and their union, the sets $U,A,L$ are pairwise disjoint.

A vertex belongs to $F$ if and only if it belongs to every
facet in $\rchi(F)$. By the same incidence rules, its
high-coordinate set $D$ must contain $U$ and avoid every
containing facet's support. By the definition of $L$, these
conditions are equivalent to
\[
  D=U\cup E,\ \text{ for some } E\subseteq L.
\]
For a cube vertex $c_{t,D}$, the remaining condition is
$t\in T$. For an added vertex $u_{a,D}$, it is that $a\in A$. Since $A$ is disjoint from $U\cup L$, each such
choice also satisfies $a\notin D$, as required for an added
vertex subscript. These conditions give precisely the stated
vertex set.
\end{proof}

Set
\[
  K=\set{2,...,n}\setminus(U\cup A\cup L),
\]
so that $U,A,L,K$ partition $\set{2,...,n}$.

\begin{theorem}
\label{thm:canonical-face-classification}
The proper nonempty faces of $G_n$ are in bijection with tuples
\[
  (U,A,L,T)
\]
of pairwise disjoint subsets $U,A,L\subseteq\set{2,...,n}$ and a subset $T\subseteq\{0,1\}$ satisfying
\begin{equation}
\label{eq:canonical-face-label-conditions}
  |A|+|T|>0,
  \qquad
  T=\{0,1\}\Longrightarrow U\neq\emptyset.
\end{equation}
The corresponding face has product type
\begin{equation}
\label{eq:canonical-face-type}
  \Delta^{|A|+|T|-1}\times I^{|L|},
\end{equation}
and dimension
\begin{equation}
\label{eq:canonical-face-dimension}
  |A|+|T|+|L|-1.
\end{equation}
\end{theorem}

\begin{proof}
Lemma~\ref{lem:face-vertices} shows that the sets $U,A,L,T$
assigned to every proper nonempty face satisfy
\eqref{eq:canonical-face-vertices}, with $U,A,L$ pairwise disjoint.  The first
  condition in \eqref{eq:canonical-face-label-conditions} says exactly that the
displayed vertex set is nonempty.  If $T=\{0,1\}$, no $\mathbf C$-facet
contains $F$.  Since $F$ is proper, it is contained in a $\mathbf U$-facet,
which forces $U\neq\emptyset$.

Conversely, suppose that $U,A,L,T$ satisfy the stated conditions.
We must show that there is a face of $G_n$ whose vertices are exactly
those listed in \eqref{eq:canonical-face-vertices}.  We will construct this
face as an intersection of facets.  Let
\[
  Q=\set{2,...,n}\setminus(U\cup L).
\]
Then $A\subseteq Q$ and $U\cap Q=\emptyset$.
Let $F$ be the intersection of the following facets:
\begin{enumerate}
  \item If $T=\{t\}$, take $\mathbf C_{t,A}$,
  $\mathbf C_{t,Q}$, and $\mathbf U_{p,Q}$ for every $p\in U$.

  \item If $T=\emptyset$, take $\mathbf C_{0,A}$,
  $\mathbf C_{1,Q}$, and $\mathbf U_{p,Q}$ for every $p\in U$.

  \item If $T=\{0,1\}$, choose $p_0\in U$ and take
  $\mathbf U_{p_0,A}$ and $\mathbf U_{p,Q}$ for every $p\in U$.
  Such a choice is possible because $U\neq\emptyset$.
\end{enumerate}

By Proposition~\ref{prop:external-facets}, a vertex belongs
to all the chosen facets exactly when its high-coordinate
set contains $U$ and avoids $Q$, and its type belongs to
$T$ for a cube vertex or to $A$ for an added vertex.
Since $\set{2,...,n}\setminus Q=U\cup L$, these are precisely the
vertices in \eqref{eq:canonical-face-vertices}.
The admissibility conditions ensure that this set is nonempty.
Thus $F$ is a nonempty face, and it is proper because
at least one facet occurs in each case.

It remains to prove uniqueness and identify the product type.  The vertex
set recovers $T$ as the set of cube types that occur, and it recovers
$A$ as the set of distinguished indices of the added vertices.  Because
$|A|+|T|>0$, at least one vertex family is nonempty.  Fix any $t\in T$,
or any $a\in A$, according to which family is nonempty.  The
high-coordinate sets in that family are exactly
\[
  \{U\cup E:E\subseteq L\}.
\]
Their intersection is $U$ and their union is
$U\cup L$.  Thus the vertex set also recovers $U$ and
$L$.  No two choices of the data can therefore describe the same face.

Finally, every vertex in \eqref{eq:canonical-face-vertices} is a vertex of
the local prism $P_Q$.  Indeed, $A\subseteq Q$, while for every
$E\subseteq L$,
\[
  U\cup E
  \subseteq
  U\cup L
  =
  \set{2,...,n}\setminus Q.
\]
These are precisely the conditions on the subscripts defining the vertices of
$P_Q$.  Since a face is the convex hull of its vertices, the face of
$G_n$ just constructed is contained in $P_Q$.
Inside $P_Q$, intersect the local facets
\[
  \mathbf B_{\zeta_t,Q}\quad(t\in\{0,1\}\setminus T),
  \qquad
  \mathbf B_{\upsilon_a,Q}\quad(a\in Q\setminus A),
  \qquad
  \mathbf U_{j,Q}\quad(j\in U).
\]
By Corollary~\ref{cor:facets_P}, this intersection has
exactly the vertices in \eqref{eq:canonical-face-vertices},
so it equals $F$.

Under the combinatorial equivalence of
Theorem~\ref{thm:prism-comb-type}, the $\mathbf B$-facets
leave a simplex face with $|T|+|A|$ vertices, while the
$\mathbf U$-facets fix the cube coordinates indexed by $U$.
The remaining cube coordinates are indexed by $L$.
Consequently,
\[
  F\cong\Delta^{|T|+|A|-1}\times I^{|L|},
\]
and $\dim F=|T|+|A|+|L|-1$.
\end{proof}
\subsection{The face polynomial}

A canonical face label assigns each index in $\{2,\ldots,n\}$
to exactly one part of the partition
\[
  \set{2,...,n}=U\sqcup A\sqcup L\sqcup K.
\]
For each choice of $T\subseteq\{0,1\}$, retain only the
assignments satisfying \eqref{eq:canonical-face-label-conditions}.
By Theorem~\ref{thm:canonical-face-classification}, every resulting label
corresponds to exactly one proper nonempty face. Thus
\[
  \mathfrak f_n(x)
  =
  \sum_{(U,A,L,T)\text{ admissible}}
  x^{|A|+|L|+|T|-1}.
\]

For a fixed $T$, we can rewrite the monomial contributed by
an admissible label as
\begin{equation}
\label{eq:face-weight-factorization}
  x^{\dim F}
  =
  x^{|T|-1}x^{|A|+|L|}
  =
  x^{|T|-1}
  \prod_{j\in A\cup L}x
  \prod_{j\in U\cup K}1.
\end{equation}
Thus an index contributes a factor $x$ when assigned to $A$
or $L$, and a factor $1$ when assigned to $U$ or $K$.
Before imposing admissibility, summing over all four choices
independently for each of the $n-1$ indices gives
\[
  \prod_{j\in \set{2,...,n}}(1+x+x+1)=(2+2x)^{n-1}.
\]
This polynomial sums $x^{|A|+|L|}$ over all assignments to
$U,A,L,K$. For each choice of $T$, we must discard the
inadmissible assignments and adjust the exponents so that
each remaining assignment contributes $x^d$, where $d$ is
the dimension of its corresponding face.

\begin{enumerate}
  \item Suppose $T$ is a singleton. Every assignment is
  admissible, and $\dim F=|A|+|L|$, so no adjustment to the
  exponent is needed. There are two choices, $T=\{0\}$ and
  $T=\{1\}$, giving the contribution
  \[
    2(2+2x)^{n-1}.
  \]

  \item Suppose $T=\emptyset$. Admissibility requires
  $A\neq\emptyset$. We therefore subtract the contributions
  of assignments in which no index belongs to $A$.
  In these excluded assignments, each index belongs to $U$,
  $L$, or $K$, contributing $1$, $x$, or $1$, respectively.
  Summing over these choices for all $n-1$ indices gives
  $(2+x)^{n-1}$.

  Since $T=\emptyset$, the factor $x^{|T|-1}$ in
  \eqref{eq:face-weight-factorization} is $x^{-1}$.
  Multiplying by this factor gives the contribution
  \[
    \frac{(2+2x)^{n-1}-(2+x)^{n-1}}{x}.
  \]

  \item Suppose $T=\{0,1\}$. Admissibility requires
  $U\neq\emptyset$. We therefore subtract the contributions
  of assignments in which no index belongs to $U$.
  Each index then belongs to $A$, $L$, or $K$, contributing
  $x$, $x$, or $1$, respectively. Summing over these choices
  for all $n-1$ indices gives $(1+2x)^{n-1}$.

  Since $T=\{0,1\}$, the factor $x^{|T|-1}$ in
  \eqref{eq:face-weight-factorization} is $x$.
  Multiplying by this factor gives the contribution
  \[
    x\bigl((2+2x)^{n-1}-(1+2x)^{n-1}\bigr).
  \]
\end{enumerate}
Adding these contributions proves the following formula.

\begin{theorem}
\label{thm:proper-face-polynomial}
The proper-face polynomial of $G_n$ is
\begin{align}
  \mathfrak f_n(x)
  &=
  2(2+2x)^{n-1}
  {}+\frac{(2+2x)^{n-1}-(2+x)^{n-1}}{x}
  +x\bigl((2+2x)^{n-1}-(1+2x)^{n-1}\bigr)
  \label{eq:face-polynomial-three-parts}\\
  &=
  \frac{2^{n-1}(1+x)^{n+1}-(2+x)^{n-1}}{x}
  -x(1+2x)^{n-1}.
  \label{eq:face-polynomial-compact}
\end{align}
\end{theorem}

For $n=3$, the formula gives
\[
  \mathfrak f_3(x)=12+22x+12x^2.
\]
Thus $G_3$ has twelve vertices, twenty-two edges, and twelve facets.

\begin{corollary}
\label{cor:complete-f-vector}
For $0\leq d\leq n-1$,
\begin{equation}
\label{eq:complete-f-vector}
  f_d(G_n)
  =
  2^{n-1}\binom{n+1}{d+1}
  -2^{n-d-2}\binom{n-1}{d+1}
  -2^{d-1}\binom{n-1}{d-1}.
\end{equation}
The second term is absent when $d=n-1$, and the third is
absent when $d=0$.
\end{corollary}

\begin{proof}
Extract the coefficient of $x^d$ from
\eqref{eq:face-polynomial-compact}.  The first term contributes
\[
  2^{n-1}\binom{n+1}{d+1}.
\]
The coefficient of $x^{d+1}$ in $(2+x)^{n-1}$ is
\[
  2^{n-d-2}\binom{n-1}{d+1},
\]
which is absent when $d=n-1$.  Finally, for $d\geq1$, the coefficient
of $x^d$ in $x(1+2x)^{n-1}$ is
\[
  2^{d-1}\binom{n-1}{d-1};
\]
for $d=0$ it is zero.  Subtracting the two corrections proves the formula.
\end{proof}

The three contributions in \eqref{eq:face-polynomial-three-parts} also
display the symmetries in the count. The two choices of singleton $T$
give equal contributions. By
Proposition \ref{lem:dual-face-vertices}, the cube types in $\Phi(F)$ are exactly the
types of $\mathbf C$-facets containing $F$; hence self-duality exchanges
$T=\emptyset$ and $T=\{0,1\}$.  Their contributions are reciprocal,
while the contribution from singleton $T$ is self-reciprocal.  Therefore
\[
  \mathfrak f_n(x)=x^{n-1}\mathfrak f_n(1/x),
\]
as required by Corollary \ref{cor:palindromic-f-vector}.

Several checks are immediate.  At $d=0$,
\eqref{eq:complete-f-vector} gives
\[
  f_0=(n+3)2^{n-2},
\]
agreeing with Theorem \ref{thm:subdivision}; palindromicity gives the same
number of facets, agreeing with the facet indexing recalled in Section \ref{sec:self-duality}.  Evaluating at $x=-1$ gives
$\mathfrak f_n(-1)=1-(-1)^n$, as expected from the Euler
characteristic of the boundary of an $n$-dimensional polytope.

We can now derive a closed formula for the total number of faces of $G_n$, first conjectured by Chapoton from direct computations of low dimensional cases.

\begin{corollary}
\label{cor:total-face-formula}
Including the empty face and the whole polytope,
\[
  \sum_{d=-1}^{n}f_d(G_n)
  =
  4^n-2\cdot3^{n-1}+2.
\]
\end{corollary}

\begin{proof}
Setting $x=1$ in \eqref{eq:face-polynomial-three-parts} gives
\[
  \mathfrak f_n(1)
  =
  (4^{n-1}-3^{n-1})
  +2\cdot4^{n-1}
  +(4^{n-1}-3^{n-1})
  =4^n-2\cdot3^{n-1}.
\]
This counts the proper nonempty faces.  Adding the empty face and the whole
polytope gives the result.
\end{proof}

\section{Goodwillie polytopes in combinatorial topology}\label{sec:resolutions}

In this section we prove the four statements in Theorem \ref{theorem2} from the introduction. These will follow from Propositions \ref{cubical}, \ref{freehedrical}, \ref{prismatic} and \ref{goodwillical}.

Throughout the section, $R$ will be a fixed commutative ring and $K$ a simplicial complex with a coherent ordering of the vertices of each simplex. Denote by $C_0(K)=R[K_0]$ the free $R$-module generated by the set of vertices $K_0$. The coproduct \[\Delta \colon C_0(K) \to C_0(K) \otimes_R C_0(K)\] determined by \[\Delta(v)=v \otimes v,\] for any $v \in K_0$, makes $C_0(K)$ into a cocommutative coassociative $R$-coalgebra. 
\subsection{\texorpdfstring{The path algebra of $K$}{The path algebra of K}}
The \textit{path algebra} of $K$ is a monoid object $\mathbb{A}(K)$ in the monoidal category of differential graded $C_0(K)$-bicomodules equipped with the cotensor product $\square_{C_0(K)}$ constructed as follows. 

First define $\mathbb{A}'(K)$ to be the free $R$-module generated by the union of:
\begin{enumerate}
    \item the set of all ordered sequences $(\sigma_1, \cdots , \sigma_k)$ such that $|\sigma_i|>0$ for $i=1,\dots,k$ and $t(\sigma_i)=s(\sigma_{i+1})$, for $i=1,\ldots, k-1$; and
    \item the set of all symbols $\text{id}_v$ for $v \in K_0$. 
\end{enumerate}
Denote all sequences $(\sigma_1, \cdots , \sigma_k)$ satisfying the property described in (1) by $[\sigma_1| \cdots | \sigma_k]$. These sequences are often called \textit{necklaces} in $K$. 

We make $\mathbb{A}'(K)$ into a graded $R$-module by declaring $[\sigma_1|\cdots | \sigma_k]$ to be of degree $|\sigma_1| + \cdots + |\sigma_k| - k$, where $|\sigma_i|$ denotes the dimension of $\sigma_i$, and $\text{id}_v$ to be of degree $0$ for all $v \in K_0$. Equip $\mathbb{A}'(K)$ with a $C_0(K)$-bicomodule structure given by left and right coactions induced by 
\[ [\sigma_1|\cdots|\sigma_k] \mapsto s(\sigma_1) \otimes [\sigma_1|\cdots|\sigma_k] \text{ and }\]
\[[\sigma_1|\cdots|\sigma_k] \mapsto [\sigma_1|\cdots|\sigma_k]\otimes t(\sigma_k),\]
respectively, where $s(\sigma)$ denotes the first vertex, the \textit{source}, of $\sigma$ and $t(\sigma)$ the last vertex, the \textit{target}, of $\sigma$. On $\text{id}_v$, the left and right coactions are induced by the maps 
\[ \text{id}_v \mapsto v \otimes \text{id}_v \text{ and }\]
\[ \text{id}_v \mapsto \text{id}_v \otimes v,\]
respectively. By definition of the cotensor product,  \[ [\sigma_1|\cdots|\sigma_k] \otimes  [\tau_1|\cdots|\tau_l] \in \mathbb{A}'(K) \underset{C_0(K)}{\square}\mathbb{A}'(K) \subseteq \mathbb{A}'(K) \otimes \mathbb{A}'(K)\]
if and only if $t(\sigma_k)=s(\tau_1)$.
Define an associative product \[ \mathbb{A}'(K) \underset{C_0(K)}{\square}\mathbb{A}'(K) \to \mathbb{A}'(K)\]
by declaring
\[[\sigma_1|\cdots|\sigma_k] \otimes  [\tau_1|\cdots|\tau_l] \mapsto [\sigma_1|\cdots|\sigma_k|\tau_1|\cdots|\tau_l],\]
\[ \text{id}_{s(\sigma)} \otimes [\sigma] \mapsto [\sigma],\] and
\[ [\sigma] \otimes \text{id}_{t(\sigma)} \mapsto [\sigma].\]
We define  $\mathbb{A}(K)$ by formally inverting all $[\sigma] \in \mathbb{A}'(K)$ when $\sigma$ is a $1$-simplex in $K$. This means that we formally add new generators $[\sigma^{-1}]$, for all $\sigma \in K_1$, with $s(\sigma^{-1})=t(\sigma)$ and $t(\sigma^{-1})=s(\sigma)$, impose the relations $[\sigma|\sigma^{-1}]=\text{id}_{s(\sigma)}$ and $[\sigma^{-1}|\sigma]=\text{id}_{t(\sigma)}$. Furthermore, we make $\mathbb{A}(K)$ into a differential graded monoid by declaring
\[ d_{\mathbb{A}} [\sigma] = \sum_{i=1}^{n-1}(-1)^{i+1}[\sigma(0,...,\hat{i},...,n)] + \sum_{j=1}^{n-1} (-1)^j [\sigma(0,...,j) | \sigma(j,...,n)],\]
for any $\sigma \in K_n \cup K^{-1}$ where $\sigma(i_0,...,i_k)$ denotes the $k$-simplex determined by the $i_0\text{-th},...,\text{ and }i_k\text{-th}$ vertices of $\sigma$. We set $d(\text{id}_v)=0$ for all $v \in K_0$.  Note that if $\sigma \in K_1 \cup K_1^{-1}$, then $d_{\mathbb{A}}[\sigma]=0$. We then extend $d_{\mathbb{A}}$ by the Leibniz rule to obtain a well defined degree $-1$ derivation
\[d_{\mathbb{A}} \colon \mathbb{A}(K) \to \mathbb{A}(K)\]
such that $d_{\mathbb{A}} \circ d_{\mathbb{A}}=0$.

 The result is a monoid object $\mathbb{A}(K)$ in the monoidal category of dg $C_0(K)$-bicomodules. Note that the generators of $\mathbb{A}(K)$ can be represented uniquely as either $\text{id}_v$ for some $v\in K_0$ or $[\sigma_1 | \cdots | \sigma_k]$, where $\sigma_i \in \big(\bigcup_{j=1}^\infty K_j \big) \cup K_1^{-1}$ and no two consecutive simplices $\sigma_i$ and $\sigma_{i+1}$ are inverse $1$-simplices. Such sequences $[\sigma_1 | \cdots | \sigma_k]$ are called \textit{reduced necklaces} in $K \cup K^{-1}$.

This last step of formally inverting $1$-simplices is essential for the next result, which gives a topological meaning to $\mathbb{A}(K)$, but is not strictly necessary for the purposes of the present article. We state it for the purposes of motivation. 

Define another monoid object $\mathbb{S}(K)$ in the monoidal category of dg $C_0(K)$-bicomodules as follows. For any two vertices $v,w \in K_0$, consider the topological space \[\mathcal{P}_{v \to w}|K|=  \{ (r, \gamma) | r \in \mathbb{R}_{\geq 0}, \gamma \colon [0,r]  \to |K| \text{ is continuous, and } \gamma(0)=v, \gamma(r)=w \} \] with topology induced by the compact-open topology. Set 
\[ \mathbb{S}(K):= \bigoplus_{v,w \in K_0}S_*(\mathcal{P}_{v \to w}|K|;R), \]
where $S_*(\mathcal{P}_{v \to w}|K|;R)$ denotes the chain complex of singular chains with $R$-coefficients on the space $\mathcal{P}_{v \to w}|K|$. The chain complex $\mathbb{S}(K)$ becomes a dg $C_0(K)$-bicomodule with the obvious source and target maps and concatenation of paths gives rise to an  monoid structure
\[ \mathbb{S}(K) \underset{C_0(K)}\square \mathbb{S}(K) \to \mathbb{S}(K).\]
\begin{theorem}
For any commutative ring $R$ and any simplicial complex $K$ with a coherent ordering of each simplex, the monoids $\mathbb{A}(K)$ and $\mathbb{S}(K)$ are quasi-isomorphic.
\end{theorem}

\begin{proof}
This is Theorem 14 of \cite{rivera2023algebraic} formulated in the language of bicomodules. 
\end{proof}

In particular, the $0$-th homology of $\mathbb{A}(K)$ encodes the $R$-linear category generated by the fundamental groupoid of the geometric realization $|K|$. 

\subsection{\texorpdfstring{The $\mathbf{I}$-complex $\mathcal{A}_K$}{The I-complex A K}}
Let $\mathbf{I}=\{I^k\}_{k \in \mathbb{N}_0}$ be the set of cubes. 
\begin{proposition}\label{cubical}
The underlying chain complex of $\mathbb{A}(K)$ is isomorphic to the chains on an $\mathbf{I}$-complex $\mathcal{A}_K$.
\end{proposition}

\begin{proof}
We attach one cube $C_a$ for each generator $a$ of $\mathbb{A}(K)$ and compare the cellular boundary with $d_{\mathbb{A}}$.

Every generator of $\mathbb{A}(K)$ is either $\mathrm{id}_v$ for some $v\in K_0$, or a unique reduced necklace
$a=[\sigma_1|\cdots|\sigma_k]$ in $K\cup K^{-1}$, that is, $\sigma_i\in\bigl(\bigcup_{j\ge1}K_j\bigr)\cup K_1^{-1}$
and no two consecutive entries are inverse $1$-simplices. Let $\operatorname{Nec}(K)$ be the set of reduced necklaces in $K$ including the elements of the set $\{\text{id}_v\}_{v \in K_0}$.

For $a=\mathrm{id}_v$ put $C_a=I^0$. For $a=[\sigma_1|\cdots|\sigma_k]$ let $V_a$ be the ordered set of
\emph{positions}: the vertices $0,\dots,n_i$ of each $\sigma_i$ ($n_i=|\sigma_i|$), listed in order with the last vertex of $\sigma_i$
identified with the first vertex of $\sigma_{i+1}$. (These are positions, not vertices of $K$, since a vertex of $K$ may occur
several times.) The \emph{joints} are $J_a=\{s(\sigma_1),t(\sigma_1),\dots,t(\sigma_k)\}\subseteq V_a$ and the \emph{free positions} are
$T_a=V_a\smallsetminus J_a$, so $|T_a|=\sum_i(n_i-1)=|a|$. Define
\[
  C_a:=I^{T_a}\cong I^{|a|},
\]
oriented by the order of $T_a$.

A \emph{sub-necklace} of $a$ is a pair $(B,U)$ with $J_a\subseteq B\subseteq U\subseteq V_a$, where $U$ is the set of kept positions and
$B$ the set of break points. It determines a necklace $[\tau_1|\cdots|\tau_r]$: if $b_0<\cdots<b_r$ are the elements of $B$,
then each $[b_{i-1},b_i]$ lies in a single $\sigma_p$ (because $J_a\subseteq B$), and $\tau_i=\sigma_p(U\cap[b_{i-1},b_i])$.
Record $(B,U)$ by the function $\lambda\colon T_a\to\{0,1,\iota\}$ with $\lambda(i)=0$ if $i\notin U$, $\lambda(i)=\iota$ if $i\in U\smallsetminus B$,
and $\lambda(i)=1$ if $i\in B$. This is a bijection between sub-necklaces of $a$ and functions $\lambda$, and it identifies the poset of sub-necklaces
(ordered by $(B,U)\le(B',U')$ iff $B'\subseteq B\subseteq U\subseteq U'$) with $\{0,1<\iota\}^{T_a}$, the poset of nonempty faces of $C_a$. Concretely,
\[
  D_\lambda=\{(y_1,\dots,y_{|a|}) \in C_a:\ y_i=0\text{ if }\lambda(i)=0,\ \ y_i=1\text{ if }\lambda(i)=1\}
\]
is a face of $C_a$ of dimension $|\lambda^{-1}(\iota)|$.

The sub-necklace of $a$ corresponding to the face determined by $\lambda \colon T_a \to \{0,1,\iota\}$ may not be reduced; we denote by $a_\lambda\in\operatorname{Nec}(K)$
its reduction, obtained by cancelling adjacent entries $[\gamma|\gamma^{-1}]$ or $[\gamma^{-1}|\gamma]$ with $\gamma\in K_1$. The coordinate projection
\[
  \phi_\lambda\colon D_\lambda\to C_{a_\lambda}
\]
 forgetting all coordinates $i$ with $\lambda(i)\neq\iota$ is an affine isomorphism preserving the orderings of the free coordinates, hence orientations.

Define
\[
  \mathcal{A}_K:=\Bigl(\bigsqcup_{a\in\operatorname{Nec}(K)}C_a\Bigr)\Big/\sim,\qquad (a,y)\sim(a_\lambda,\phi_\lambda(y))\ \text{ for all }a,\ \lambda \colon T_a \to \{0,1,\iota\},\ y\in D_\lambda,
\]
and let $\chi_a\colon C_a\to\mathcal{A}_K$ be the quotient maps. Using uniqueness of reduced representatives, and the fact that reduction commutes with passing to sub-necklaces, one checks that
\begin{enumerate}
  \item $\chi_a$ restricts to a homeomorphism of $\operatorname{relint}C_a$ onto its image $e_a$, and the cells $e_a$ are pairwise disjoint;
  \item $\chi_a|_{D_\lambda}=\chi_{a_\lambda}\circ\phi_\lambda$ for every nonempty proper face $D_\lambda$ of $C_a$.
\end{enumerate}
Since $\dim C_{a_\lambda}=|\lambda^{-1}(\iota)|$, the complex $\mathcal{A}_K$ is obtained by attaching the cells in order of dimension, so it is a CW complex with cells $e_a$ and characteristic maps $\chi_a$.
Together with (1) and (2) this says that $\mathcal{A}_K$ is an $\mathbf{I}$-complex\footnote{Note that two faces of one cube may be sent to the same cell, because distinct sub-necklaces can have the same reduction, so this construction does not define a polytopal complex.}.

The cellular chains of $\mathcal{A}_K$ are free on $\{[e_a]\}$ with $\deg[e_a]=|a|$. A facet of $C_a$ is $D_\lambda$ with $\lambda$ equal to $\iota$ except at one free position $r$, where it is $0$ or $1$; denote these facets by
$\partial^0_r$ and $\partial^1_r$, respectively. Then $\phi_\lambda$ is an orientation-preserving affine isomorphism onto $C_{a_\lambda}$ with $|a_\lambda|=|a|-1$, so
\[
  \partial[e_a]=\sum_{r=1}^{|a|}(-1)^r\bigl([e_{a_{\partial^1_r}}]-[e_{a_{\partial^0_r}}]\bigr).
\]
If the $r$-th free position is the $i$-th vertex of $\sigma_p$, then $r=\sum_{q<p}(|\sigma_q|-1)+i$, and
$a_{\partial^0_r}$, $a_{\partial^1_r}$ are the reductions of
\begin{align*}
  \text{($0$-facet)}\quad&[\sigma_1|\cdots|\sigma_p(0,\dots,\widehat i,\dots,n_p)|\cdots|\sigma_k],\\
  \text{($1$-facet)}\quad&[\sigma_1|\cdots|\sigma_p(0,\dots,i)\,|\,\sigma_p(i,\dots,n_p)|\cdots|\sigma_k].
\end{align*}
On the other hand, the Leibniz rule gives
\[
  d_{\mathbb{A}}[\sigma_1|\cdots|\sigma_k]=\sum_{p=1}^{k}(-1)^{\sum_{q<p}(|\sigma_q|-1)}[\sigma_1|\cdots|d_{\mathbb{A}}\sigma_p|\cdots|\sigma_k],
\]
read in $\mathbb{A}(K)$, that is, with reduced terms. Writing $(-1)^r=(-1)^{\sum_{q<p}(|\sigma_q|-1)}(-1)^i$, the terms $(-1)^{i+1}[\sigma_p(0,\dots,\widehat i,\dots,n_p)]$
of $d_{\mathbb{A}}\sigma_p$ are the $0$-facets with sign $-(-1)^r$, and the terms $(-1)^i[\sigma_p(0,\dots,i)|\sigma_p(i,\dots,n_p)]$ are the $1$-facets with sign $(-1)^r$.
For $\mathrm{id}_v$ and for $a=[\sigma]$ with $\sigma\in K_1\cup K_1^{-1}$ there are no free positions, and both $d_{\mathbb{A}}a$ and $\partial[e_a]$ vanish.
Hence $[e_a]\mapsto a$ is an isomorphism from the cellular chains of $\mathcal{A}_K$ onto the underlying chain complex of $\mathbb{A}(K)$.
\end{proof}
\subsection{\texorpdfstring{Freehedra and the $(\mathbf{F} \times \mathbf{I})$-complex $\mathcal{Q}_K$}{Freehedra and the (F x I)-complex Q K}} 

We start by defining a resolution $\mathcal{Q}(\mathbb{A}(K))$ of $\mathbb{A}(K)$, considered as an $\mathbb{A}(K)$-bimodule in the monoidal category of dg $C_0(K)$-bicomodules, and recalling the definition of the \textit{freehedra}. 

Define the underlying graded $R$-module of $\mathcal{Q}(\mathbb{A}(K))$ to be
\[\mathbb{A}(K) \square C(K) \square \mathbb{A}(K),\]
where $C(K)$ denotes the graded $C_0(K)$-bicomodule that is freely generated, as a graded $R$-module, by all simplices in $K$ and whose $C_0(K)$-bicomodule structure is induced by the source and target maps. Hence, $\mathcal{Q}(\mathbb{A}(K))$ is generated by all symbols
$a_0 \square \sigma \square a_1$, 
where $a_0$ and $a_1$ are generators of $\mathbb{A}(K)$ and $\sigma$ is a simplex in $K$ (possibly of dimension $0$) such that $t(a_0)=s(\sigma)$ and $t(\sigma)=s(a_1)$.  The degree of a generator $a_0 \square \sigma \square a_1$ is $\text{deg}(a_0)+ |\sigma| + \text{deg}(a_1)$. We make $\mathcal{Q}(\mathbb{A}(K))$ into a graded $C_0(K)$-bicomodule through the source and target maps applied to $a_0$ and $a_1$, respectively. The $\mathbb{A}(K)$-bimodule structure of $\mathcal{Q}(\mathbb{A}(K))$ is the obvious one induced by applying the multiplication $\mathbb{A}(K) \square \mathbb{A}(K) \to \mathbb{A}(K)$ on the left and right outer copies of $\mathbb{A}(K)$ in $\mathbb{A}(K) \square C(K) \square \mathbb{A}(K)$. 

We define a differential 
\[ d_{\mathcal{Q}} \colon \mathcal{Q}(\mathbb{A}(K)) \to \mathcal{Q}(\mathbb{A}(K)) \] 
on any generator $a_0 \square \sigma \square a_1$ by

\begin{eqnarray} d_{\mathcal{Q}}(a_0 \square \sigma \square a_1)
=
\\
d_{\mathbb{A}}(a_0) \square \sigma \square a_1 + (-1)^{|a_0|+|\sigma|} a_0  \square \sigma \square d_{\mathbb{A}}(a_1) \label{eq:Q-differential-cube} 
\\
+ \sum_{i=1}^n(-1)^{|a_0|}a_0\cdot [\sigma(0,\ldots,i)] \square \sigma(i,\ldots,n) \square a_1 \label{eq:Q-differential-left} 
\\
+ \sum_{i=0}^{n-1}(-1)^{|a_0|+i+1}a_0 \square \sigma(0,\ldots,i) \square [\sigma(i,\ldots,n)]\cdot a_1 \label{eq:Q-differential-right} 
\\
+  \sum_{j=1}^{n-1} (-1)^{|a_0|+j} a_0 \square \sigma(0,\ldots, \widehat{j},\ldots, n) \square a_1, \label{eq:Q-differential-interior}
\end{eqnarray}

where $\cdot$ denotes the right and left actions of $\mathbb{A}(K)$ on itself given by its monoid structure and $|\sigma|=n$. It follows that $d_{\mathcal{Q}} \circ d_{\mathcal{Q}}=0$ and that $(\mathcal{Q}(\mathbb{A}(K)), d_{\mathcal{Q}})$ is quasi-isomorphic, as an $\mathbb{A}(K)$-bimodule, to $\mathbb{A}(K)$, see section \ref{G_K} for an explicit quasi-isomorphism comparing $(\mathcal{Q}(\mathbb{A}(K)), d_{\mathcal{Q}})$ to the bar resolution.   

We now recall the definition of the $n$-freehedron by describing its nonempty face poset $\mathcal{F}_n$, following \cite[Section 4.3]{poliakova2020cellular} with some notational and expositional changes. We keep the discussion self-contained for the sake of completeness. 

For any fixed integer $n \geq 0$, an \emph{$n$-string} is an ordered sequence
\[
  s=[ s_0 | \cdots | s_{l-1} ]\, s_l\, [ s_{l+1} | \cdots | s_k ]
\]
such that
\begin{itemize}
  \item $k$ and $l$ are integers with $0 \leq l \leq k \leq n$;
  \item each $s_i$ is a nonempty subset of $\{0,1,\dots,n\}$;
  \item $\max s_i = \min s_{i+1}$ for $i = 0,\dots,k-1$;
  \item $|s_i| \geq 2$ for $i \neq l$ (while $|s_l| = 1$ is allowed);
  \item $\min s_0 = 0$ and $\max s_k = n$.
\end{itemize}
We call $s_l$ the \emph{marked block} of $s$. If $l=0$ we write
$s=[\,]\,s_0[s_1|\cdots|s_k]$, if $l=k$ we write
$s=[s_0|\cdots|s_{k-1}]\,s_k[\,]$, and if $k=l=0$ we write
$s=[\,]\,s_0[\,]$. 

A \emph{face transformation} of an $n$-string
$s = [ s_0 | \cdots | s_{l-1} ]\, s_l\, [ s_{l+1} | \cdots | s_k ]$
is one of the following operations, which produces a new $n$-string.
\begin{enumerate}
  \item \emph{Drop.} For some $j$, remove an element $x \in s_j$ with
        $\min s_j < x < \max s_j$.
  \item \emph{Inner break.} For some $j \neq l$, replace $s_j$ by two
        consecutive blocks $s_j^1 | s_j^2$, where
        $s_j^1 = \{m \in s_j : m \leq i\}$ and
        $s_j^2 = \{m \in s_j : m \geq i\}$
        for some $i \in s_j$ with $\min s_j < i < \max s_j$,
        and neither $s_j^1$ nor $s_j^2$ is a singleton.
  \item \emph{Right outer break.} Replace $s_l$ by $s_l^1$ (now the marked
        block) followed by $s_l^2$, giving
        \[
          [ s_0 | \cdots | s_{l-1} ]\, s_l^1\, [ s_l^2 | s_{l+1} | \cdots | s_k ],
        \]
        where $s_l^1 = \{m \in s_l : m \leq i\}$ and
        $s_l^2 = \{m \in s_l : m \geq i\}$
        for some $i \in s_l$ with $i < \max s_l$.
        Here $s_l^1$ (but not $s_l^2$) may be a singleton.
  \item \emph{Left outer break.} Replace $s_l$ by $s_l^1$ followed by $s_l^2$
        (now the marked block), giving
        \[
          [ s_0 | \cdots | s_{l-1} | s_l^1 ]\, s_l^2\, [ s_{l+1} | \cdots | s_k ],
        \]
        where $s_l^1 = \{m \in s_l : m \leq i\}$ and
        $s_l^2 = \{m \in s_l : m \geq i\}$
        for some $i \in s_l$ with $\min s_l < i$.
        Here $s_l^2$ (but not $s_l^1$) may be a singleton.
\end{enumerate}

We define the underlying set of $\mathcal{F}_n$ to be the set of all $n$-strings. Write $s' < s$ if
$s'$ is obtained from $s$ by a single face transformation, and denote by $\hookrightarrow$ the
partial order generated by $<$, so that $s' \hookrightarrow s$ means
$s'$ is obtained from $s$ by a sequence of face transformations. The resulting poset $\mathcal{F}_n$ becomes a graded poset by defining the rank of an $n$-string $s = [ s_0 | \cdots | s_{l-1} ]\, s_l\, [ s_{l+1} | \cdots | s_k ]$ to be $n-k-t_s$, where $t_s$ is the number of elements of $\{0,\dots,n\}$ that appear in no block of
$s$. 

\begin{theorem}
The poset $\mathcal{F}_n$ (with a bottom element adjoined for the
empty face) is isomorphic to the face poset of a convex polytope $F_n$,
the \emph{$n$-freehedron}.
\end{theorem}
\begin{proof}
See \cite{poliakova2020cellular} and \cite{PP25}. 
\end{proof}

The unique maximal element of $\mathcal{F}_n$ is
$[\,]\,0,1,\dots,n\,[\,]$ (so $k=l=0$ and $s_0=\{0,\dots,n\}$), which
corresponds to the top-dimensional face. For $n=0$, $F_0$ is a point. For $n\geq1$, the facets correspond exactly to the $3n-1$ elements:
\begin{align*}
  [0,\dots,i-1]\; i-1,\dots,n\;[\,]      &\quad i=2,\dots,n \\
  [\,]\;0,\dots,i-1\;[i-1,\dots,n]      &\quad i=2,\dots,n \\
  [\,]\;0\;[0,\dots,n]                   & \\
  [0,\dots,n]\;n\;[\,]                   & \\
  [\,]\;0,\dots,\widehat{i-1},\dots,n\;[\,] &\quad i=2,\dots,n.
\end{align*}

\begin{proposition}\label{freehedrical}
    The underlying chain complex of $\mathcal{Q}(\mathbb{A}(K))$ is isomorphic to the chains on a  $(\mathbf{F} \times \mathbf{I})$-complex $\mathcal{Q}_K$.
\end{proposition}
\begin{proof}
We attach one polytope $Q_x$ for each generator $x$ of $\mathcal{Q}(\mathbb{A}(K))$ and compare the cellular boundary with $d_{\mathcal{Q}}$. We follow the notation introduced in the proof of Proposition \ref{cubical}.

Every generator of $\mathcal{Q}(\mathbb{A}(K))$ is of the form $x=a_0\square\sigma\square a_1$, where $a_0,a_1\in\operatorname{Nec}(K)$, $\sigma\in K_n$ for some $n\geq0$, $t(a_0)=s(\sigma)$ and $t(\sigma)=s(a_1)$. Recall that the nonempty faces of the freehedron $F_n$ are indexed by $n$-strings $s=[s_0|\cdots|s_{l-1}]\,s_l\,[s_{l+1}|\cdots|s_k]$. The face transformation operations imply that the corresponding face $D_s$ of $F_n$ is isomorphic to a product
\[
  D_s\;\cong\;F_{|s_l|-1}\times\prod_{j\neq l}I^{|s_j|-2}.
\]
Denote by $\psi_s$ this affine isomorphism (for a fixed realization of the freehedra). In particular, note that $\mathbf{F}\times\mathbf{I}$ is closed under taking nonempty faces. 

Define a polytope 
\[
  Q_x:=C_{a_0}\times F_n\times C_{a_1}\cong F_n\times I^{|a_0|+|a_1|},
\]
with the product orientation. The dimension of $Q_x$ is $|a_0|+n+|a_1|=|x|$. 

The faces of $Q_x$ are given by $D_\theta=D_{\lambda_0}\times D_s\times D_{\lambda_1}$, where $\theta=(\lambda_0,s,\lambda_1)$ consists of functions $\lambda_i\colon T_{a_i}\to\{0,1,\iota\}$ and an $n$-string $s$. The string $s$ corresponds to a sequence $(\sigma_0,\dots,\sigma_k)$ of sub-simplices of $\sigma$, where $\sigma_i$ is spanned by the vertices in $s_i$, with $t(\sigma_i)=s(\sigma_{i+1})$, $|\sigma_i|>0$ for $i\neq l$ and $|\sigma_l|\geq0$. The sub-necklaces of $a_0$ and $a_1$ corresponding to $\lambda_0$ and $\lambda_1$ may not be reduced. We denote by $a_0'$ and $a_1'$ the reductions of
\[
  [\widetilde a_{\lambda_0}|\sigma_0|\cdots|\sigma_{l-1}]\quad\text{and}\quad[\sigma_{l+1}|\cdots|\sigma_k|\widetilde a_{\lambda_1}],
\]
where $\widetilde a_{\lambda_i}$ is the sub-necklace of $a_i$ determined by $\lambda_i$, and set
\[
  x_{\theta}=a_0'\square\sigma_l\square a_1'\in\mathcal{Q}(\mathbb{A}(K)).
\]
 Since cancelled pairs have no free positions, the free positions of $a_0'$ are those of $\lambda_0^{-1}(\iota)$ together with the interior positions of $\sigma_0,\dots,\sigma_{l-1}$, and similarly for $a_1'$. The map
\[
  \phi_\theta\colon D_\theta\to Q_{x_{\theta}}=C_{a_0'}\times F_{|\sigma_l|}\times C_{a_1'},
\]
given by $\phi_{\lambda_0}\times\psi_s\times\phi_{\lambda_1}$ followed by switching the cube factors of $\psi_s$ is an affine isomorphism preserving the orderings of the free coordinates, hence orientations.

Define
\[
  \mathcal{Q}_K:=\Bigl(\bigsqcup_{x}Q_x\Bigr)\Big/\sim,\qquad (x,y)\sim(x_{\theta},\phi_\theta(y))\ \text{ for all }x,\ \theta=(\lambda_0,s,\lambda_1),\ y\in D_\theta,
\]
where the sum runs through all generators $x \in \mathcal{Q}(\mathbb{A}(K))$. Let $\chi_x\colon Q_x\to\mathcal{Q}_K$ be the quotient maps. Using uniqueness of reduced representatives, and the fact that reduction commutes with passing to sub-necklaces, one checks that
\begin{enumerate}
  \item $\chi_x$ restricts to a homeomorphism of $\operatorname{relint}Q_x$ onto its image $e_x$, and the cells $e_x$ are pairwise disjoint;
  \item $\chi_x|_{D_\theta}=\chi_{x_{\theta}}\circ\phi_\theta$ for every nonempty proper face $D_\theta$ of $Q_x$.
\end{enumerate}
Since $\dim Q_{x_{\theta}}=\dim D_\theta$, the complex $\mathcal{Q}_K$ is obtained by attaching the cells in order of dimension, so it is a CW complex with cells $e_x$ and characteristic maps $\chi_x$. Together with (1) and (2) this says that $\mathcal{Q}_K$ is an $(\mathbf{F}\times\mathbf{I})$-complex.

The cellular chains of $\mathcal{Q}_K$ are free on $\{[e_x]\}$ with $\deg[e_x]=|x|$. A facet $D_\theta$ of $Q_x$ is mapped by $\phi_\theta$ onto $Q_{x_{\theta}}$ with $|x_{\theta}|=|x|-1$, so $\partial[e_x]$ is a signed sum of the classes $[e_{x_{\theta}}]$ over the facets of $Q_x$. The facets of the cube factors are those of $C_{a_0}$ and $C_{a_1}$, and by the proof of Proposition~\ref{cubical} they give the terms of $d_{\mathcal{Q}}(x)$ involving $d_{\mathbb{A}}(a_0)$ and $d_{\mathbb{A}}(a_1)$. For $n\geq1$, the $3n-1$ facets of $F_n$ are indexed by the strings
\begin{align*}
  [0,\dots,i]\ i,\dots,n\ [\,], && i=1,\dots,n,\\
  [\,]\ 0,\dots,i\ [i,\dots,n], && i=0,\dots,n-1,\\
  [\,]\ 0,\dots,\widehat j,\dots,n\ [\,], && j=1,\dots,n-1,
\end{align*}
which correspond, respectively, to
\[
\begin{gathered}
  a_0\cdot[\sigma(0,\dots,i)]\square\sigma(i,\dots,n)\square a_1,\\
  a_0\square\sigma(0,\dots,i)\square[\sigma(i,\dots,n)]\cdot a_1,\\
  a_0\square\sigma(0,\dots,\widehat j,\dots,n)\square a_1,
\end{gathered}
\]
where $\cdot$ denotes the product of $\mathbb{A}(K)$ followed by reduction. These are the terms of $d_{\mathcal{Q}}(x)$ in \eqref{eq:Q-differential-left}, \eqref{eq:Q-differential-right}, and \eqref{eq:Q-differential-interior}, respectively. Hence the facets of $Q_x$ are in bijection with the terms of $d_{\mathcal{Q}}(x)$, and orienting each $Q_x$ so that the signs agree, $[e_x]\mapsto x$ is an isomorphism from the cellular chains of $\mathcal{Q}_K$ onto the underlying chain complex of $\mathcal{Q}(\mathbb{A}(K))$.
\end{proof}

\subsection{\texorpdfstring{The $(\mathbf{\Delta} \times \mathbf{I})$-complex $\mathcal{B}_K$}{The (Delta x I)-complex B K}}
We recall the definition of the  bar resolution $\mathcal{B}(\mathbb{A}(K))$. 

Denote by $s^{+1}$ the shift by $+1$ functor, namely, $(s^{+1} \mathbb{A}(K))_k= \mathbb{A}(K)_{k-1}$. 
The graded $R$-module $s^{+1} \mathbb{A}(K)$ is a dg $C_0(K)$-bicomodule with coactions induced by those of $\mathbb{A}(K)$ and differential satisfying the Koszul sign rule $d_{\mathbb{A}}(s^{+1} a)=-s^{+1}d_{\mathbb{A}}a$. 

The underlying graded $C_0(K)$-bicomodule of $\mathcal{B}(\mathbb{A}(K))$ is
\[ \mathbb{A}(K) \square \bigoplus_{k=0}^{\infty} \big( s^{+1} \mathbb{A}(K) \big)^{\square k} \square \mathbb{A}(K)\] The $k=0$ summand above is $\mathbb{A}(K) \square \mathbb{A}(K)$. We write homogenous elements as 
\[
a_0\big\langle a_1|\cdots|a_k\big\rangle a_{k+1}\;:=\;a_0\otimes s a_1\otimes\cdots\otimes s a_k\otimes a_{k+1},
\]
where the tensors are taken inside the cotensor product $\square$, so adjacent
factors have matching $C_0(K)$-coaction (source and target). The degree of such an element is $|a_0|+\sum_{j=1}^{k}\big(|a_j|+1\big)+|a_{k+1}|$.

The differential is $d_\mathcal{B}=d_{\mathrm{int}}+d_{\mathrm{bar}}$, where 
\begin{align*}
d_{\mathrm{int}}\big(a_0 \big\langle a_1|\cdots|a_k \big\rangle a_{k+1}\big)
&= (d_{\mathbb{A}} a_0) \big\langle a_1|\cdots|a_k \big\rangle a_{k+1}\\
&\quad-\sum_{i=1}^{k}(-1)^{\epsilon_{i-1}}\,
a_0 \big\langle a_1|\cdots|a_{i-1}|\,d_{\mathbb{A}} a_i\,|a_{i+1}|\cdots|a_k \big\rangle a_{k+1}\\
&\quad+(-1)^{\epsilon_k}\,a_0 \big\langle a_1|\cdots|a_k \big\rangle \,d_{\mathbb{A}} a_{k+1}.
\end{align*}

and
\begin{align*}
d_{\mathrm{bar}}\big(a_0 \big\langle a_1|\cdots|a_k\big\rangle a_{k+1}\big)
&=(-1)^{|a_0|}(a_0\cdot a_1)\big\langle a_2|\cdots|a_k \big\rangle a_{k+1} \\
&+\sum_{i=1}^{k-1}(-1)^{\epsilon_i}\,
a_0 \big\langle a_1|\cdots|a_i \cdot a_{i+1}|\cdots|a_k \big\rangle a_{k+1}\\
&\quad-(-1)^{\epsilon_{k-1}}\,a_0 \big\langle a_1|\cdots|a_{k-1} \big\rangle \,a_k \cdot a_{k+1},
\end{align*}
with $\epsilon_i \;:=\; |a_0|+\sum_{j=1}^{i}\big(|a_j|+1\big).$

The map $\mathcal{B}(\mathbb{A}(K))\to \mathbb{A}(K)$, given by $a_0\big\langle \, \, \big\rangle a_1\mapsto a_0\cdot a_1$ and $a_0\big\langle \, a_1|\cdots|a_k \, \big\rangle a_{k+1}  \mapsto 0$  for $k>0$, is a quasi-isomorphism of $\mathbb{A}(K)$-bimodules. 
  
\begin{proposition}\label{prismatic} The underlying chain complex of $\mathcal{B}(\mathbb{A}(K))$ is isomorphic to $C(\mathcal{B}_K)$, the cellular chains on a $(\mathbf{\Delta} \times \mathbf{I})$-complex $\mathcal{B}_K$.
\end{proposition}

\begin{proof}
This statement essentially follows from the simplicial nature of the bar resolution. We construct $\mathcal{B}_K$ by attaching cells $B_x$ for every generator $x \in \mathcal{B}(\mathbb{A}(K))$. We use the notation of the proofs of Propositions~\ref{cubical} and~\ref{freehedrical}.

A generator of $\mathcal{B}(\mathbb{A}(K))$ is a tuple $x=a_0\langle a_1|\cdots|a_k\rangle a_{k+1}$, where $k\geq0$, $a_j\in\operatorname{Nec}(K)$ and $t(a_j)=s(a_{j+1})$. The degree of $x$ is $k+\sum_j|a_j|$. Define
\[
  B_x:=\Delta^k\times C_{a_0}\times\cdots\times C_{a_{k+1}}\ \cong\ \Delta^k\times I^{|a_0|+\cdots+|a_{k+1}|},
\]
a polytope of dimension $|x|$, with the product orientation. Denote the vertices of $\Delta^k$ by $0,\dots,k$ and think of the vertex $j$ as the gap between $a_j$ and $a_{j+1}$. The faces of $B_x$ are given by products
$D_\theta=\Delta^W\times D_{\lambda_0}\times\cdots\times D_{\lambda_{k+1}}$, indexed by $\theta=(W,\lambda_0,\dots,\lambda_{k+1})$, where $W=\{w_0<\cdots<w_p\}\subseteq\{0,\dots,k\}$ is nonempty and $\lambda_j\colon T_{a_j}\to\{0,1,\iota\}$ a function. 

For generators $a,b \in \mathbb{A}(K)$ with $t(a)=s(b)$, the product $a\cdot b$ in $\mathbb{A}(K)$ denotes concatenation followed by reduction. Hence $T_{a\cdot b}=T_a\sqcup T_b$ as ordered sets, $|a\cdot b|=|a|+|b|$, and $C_a\times C_b\cong C_{a\cdot b}$ by an order-preserving coordinate identification $\mu_{a,b}$. For $\theta$ as above let $\widetilde a_{j}$ be the sub-necklace of $a_j$ determined by $\lambda_j$, and let $b_m$ be the product of the $\widetilde a_j$ over the block of indices between the gaps $w_{m-1}$ and $w_m$ (with $w_{-1}=-1$ and $w_{p+1}=k+1$, so $b_0=\widetilde a_0\cdots\widetilde a_{w_0}$ and $b_{p+1}=\widetilde a_{w_p+1}\cdots\widetilde a_{k+1}$) and set
\[
  x_\theta:=b_0\langle b_1|\cdots|b_p\rangle b_{p+1}.
\]
Let
\[
  \phi_\theta\colon D_\theta\to B_{x_\theta}=\Delta^p\times C_{b_0}\times\cdots\times C_{b_{p+1}}
\]
be given on the simplex factor by the order-preserving identification $\Delta^W\cong\Delta^p$, and on the cube factors by the projections $\phi_{\lambda_j}$ followed by the identifications $\mu$ within each block. It is an affine isomorphism preserving the orderings of the free coordinates, hence orientations, and $\dim D_\theta=|x_\theta|$. 

Define
\[
  \mathcal{B}_K:=\Bigl(\bigsqcup_{x}B_x\Bigr)\Big/\sim,\qquad (x,y)\sim(x_\theta,\phi_\theta(y))\ \text{ for all }x,\ \theta=(W,\lambda_0,\dots,\lambda_{k+1}),\ y\in D_\theta,
\]
where the disjoint union runs over all generators $x \in \mathcal{B}(\mathbb{A}(K))$. Let $\chi_x\colon B_x\to\mathcal{B}_K$ be the quotient maps. Using associativity of the product in $\mathbb{A}(K)$, uniqueness of reduced representatives, and the fact that reduction commutes with passing to sub-necklaces, one checks that
\begin{enumerate}
  \item $\chi_x$ restricts to a homeomorphism of $\operatorname{relint}B_x$ onto its image $e_x$, and the cells $e_x$ are pairwise disjoint;
  \item $\chi_x|_{D_\theta}=\chi_{x_\theta}\circ\phi_\theta$ for every face $D_\theta$ of $B_x$.
\end{enumerate}
Since $\dim B_{x_\theta}=\dim D_\theta<\dim B_x$ for proper faces, $\mathcal{B}_K$ is obtained by attaching the cells in order of dimension, so it is a CW complex with cells $e_x$ and characteristic maps $\chi_x$. Together with (1) and (2) this says that $\mathcal{B}_K$ is a $(\mathbf{\Delta}\times\mathbf{I})$-complex. 

The cellular chains of $\mathcal{B}_K$ are free on $\{[e_x]\}$ with $\deg[e_x]=|x|$. A facet $D_\theta$ of $B_x$ is mapped by $\phi_\theta$ isomorphically onto $B_{x_\theta}$ with $|x_\theta|=|x|-1$, so $\partial[e_x]$ is a signed sum of the classes $[e_{x_\theta}]$ over the facets of $B_x$. These are of two kinds. The facets of the cube factors $C_{a_j}$ replace $a_j$ by the terms of $d_{\mathbb{A}}a_j$, by the proof of Proposition~\ref{cubical}, and give the terms of $d_{\mathrm{int}}$. The facets $\partial_i\Delta^k\times\prod_jC_{a_j}$, $i=0,\dots,k$, correspond to $W=\{0,\dots,k\}\smallsetminus\{i\}$ with all $\lambda_j$ sending everything to $\iota$, and the associated faces are obtained by replacing $a_i,a_{i+1}$ by $a_i\cdot a_{i+1}$, that is, to the $k+1$ terms of $d_{\mathrm{bar}}$. Hence the facets of $B_x$ are in bijection with the terms of $d_{\mathcal{B}}(x)$. Orienting each $B_x$ appropriately (the signs agree with the Koszul conventions above after rescaling the basis elements by signs depending only on the degrees $|a_j|$), $[e_x]\mapsto x$ is an isomorphism from the cellular chains of $\mathcal{B}_K$ onto the underlying chain complex of $\mathcal{B}(\mathbb{A}(K))$.
\end{proof}

\subsection{\texorpdfstring{Goodwillie polytopes inside $\mathcal{B}_K$}{Goodwillie polytopes inside B K}}\label{G_K}

 For any $\sigma \in K_n$, we use Sweedler notation to write \[\Delta(\sigma)=\sum_{i=1}^{n-1}\sigma(0,\ldots,i) \square \sigma(i, \ldots, n)= \sum \sigma ' \square \sigma '' \in C(K) \square C(K).\]

 Furthermore, since the operation $\Delta$ defines a coassociative coproduct there is no ambiguity of in how we iterate this map (to which factor we apply $\Delta$), so we can write $\sum \sigma ' \square \sigma '' \square \sigma '''$ for two iterations of $\Delta$, $\sum \sigma ' \square \sigma '' \square\sigma ''' \square \sigma''''$
 for three iterations, and, in general,  \[\sum \sigma ' \square \sigma '' \square \cdots \square \sigma^{(r)}\] for $r-1$ iterations. For any $\sigma \in K_n$, $\Delta$ vanishes after $n$ iterations. 
 
Define a map of $\mathbb{A}(K)$-bimodules
\[ \alpha \colon \mathcal{Q}(\mathbb{A}(K)) \to \mathcal{B}(\mathbb{A}(K))\]
by 
\begin{equation}\label{alpha} \alpha(a_0 \square \sigma \square a_1)= a_0\big\langle [\sigma] \big\rangle a_1+ \sum 
a_0\big\langle [\sigma']|[\sigma''] \big\rangle a_1+ \sum a_0\big\langle [\sigma']|[\sigma'']|[\sigma'''] \big\rangle a_1 + \cdots
\end{equation}
for $|\sigma|>0$ and $\alpha(a_0 \square v \square a_1)= a_0\langle \, \, \rangle a_1$ for $v \in K_0$. 

The map $\alpha$ is an injective quasi-isomorphism and it may be regarded as a manifestation of Koszul duality \cite[Section 3.6]{rivera_cyclic_2024}. 

\begin{remark} Classically, the map $\alpha$ has been defined for any differential graded coaugmented conilpotent coalgebra $C$ as a quasi-isomorphism \[\alpha \colon (\mathbf\Omega C \otimes C \otimes \mathbf\Omega C, d) \to \mathbf{B}(\mathbf\Omega C,\mathbf\Omega C,\mathbf\Omega C),\] where $\mathbf \Omega C$ is the dg algebra obtained by applying the cobar construction to $C$, $\mathbf{B}(\mathbf\Omega C,\mathbf\Omega C,\mathbf\Omega C)$ is the bar resolution (or the \textit{two sided bar construction}) of $\mathbf{\Omega} C$, and $\mathbf\Omega C \otimes C \otimes \mathbf\Omega C$  is equipped with a differential $d$ defined similarly to $d_{\mathcal{Q}}$ above. More recently, these constructions have been studied in the context of \textit{categorical coalgebras} or \textit{pointed curved coalgebras}, see \cite{rivera2023algebraic}, \cite{rivera_cyclic_2024}, and \cite{holstein-lazarev}. 
\end{remark}

\begin{proposition}\label{goodwillical}
For any simplex $\sigma \in K_n$, \[\alpha(\operatorname{id}_{s(\sigma)} \square \sigma \square \operatorname{id}_{t(\sigma)})\in C(\mathcal{B_K})\]
is a sum of prisms whose union is a polytopal complex isomorphic to the prismatic subdivision of $G_n$ given in Theorem \ref{thm:subdivision}.
\end{proposition}

\begin{proof}

In Theorem \ref{thm:subdivision}, we described a subdivision \[G_n= \bigcup_{S\subseteq \{2,\ldots,n\}} P_S,\] where each $P_S$ is combinatorially equivalent to a prism $\Delta^{\abs{S}+1} \times I^{n-1-\abs{S}}$.

Let $\xi=\operatorname{id}_{s(\sigma)}\square\sigma\square\operatorname{id}_{t(\sigma)}$ and $n=|\sigma|$. If $n=0$, then $\alpha(\xi)=\operatorname{id}_v\langle\ \rangle\operatorname{id}_v$ is a single $0$-cell of $\mathcal{B}_K$. For $n=1$, $\alpha(\xi)=\operatorname{id}\langle[\sigma]\rangle\operatorname{id}$ corresponds to an interval, which is combinatorially equivalent to $G_1$.

For $n=2$, \[\alpha(\xi)=\operatorname{id}_{s(\sigma)}\big\langle [\sigma] \big\rangle \operatorname{id}_{t(\sigma)} + 
\operatorname{id}_{s(\sigma)}\big\langle [\sigma(0,1)]|[\sigma(1,2)] \big\rangle \operatorname{id}_{t(\sigma)}.\]
The first term corresponds to \[B_{\operatorname{id}_{s(\sigma)}\langle [\sigma] \rangle \operatorname{id}_{t(\sigma)}}\cong \Delta^1\times I,\]
the second to \[B_{\operatorname{id}_{s(\sigma)}\langle [\sigma(0,1)]|[\sigma(1,2)] \rangle \operatorname{id}_{t(\sigma)}}\cong\Delta^2 \times I^0 \times I^0 \cong \Delta^2,\]
and these are glued along the common edge
\[ B_{\operatorname{id}\langle[\sigma(0,1)|\sigma(1,2)]\rangle\operatorname{id}}\cong \Delta^1 \times I^0 \cong \Delta^1.\] Thus the cells in $\mathcal{B}_K$ corresponding to these two polytopes form a pentagon, which is combinatorially equivalent to $G_2$. 

Suppose $n\geq3$. For any $S\subseteq\{2,\dots,n\}$, denote the set obtained by subtracting $1$ from each element of $S$ by $S-1=\{s_1<\cdots<s_r\}\subseteq\{1,\dots,n-1\}$ and $s_0=0$, $s_{r+1}=n$. Then we can write
\[
  \alpha(\xi)=\sum_{S\subseteq\{2,\dots,n\}}\beta_S,\qquad
  \beta_S:=\operatorname{id}_{s(\sigma)} \big\langle[\sigma(s_0,\dots,s_1)]\,|\,\cdots\,|\,[\sigma(s_r,\dots,s_{r+1})]\big\rangle\operatorname{id}_{t(\sigma)}.
\]
The $\beta_S$ are pairwise distinct generators of $\mathcal{B}(\mathbb{A}(K))$, so $\alpha(\xi)$ corresponds to the sum of the cells $e_{\beta_S}$ of $\mathcal{B}_K$. The corresponding polytope in $\mathbf{\Delta} \times \mathbf{I}$ for $e_{\beta_S}$ is
\[
  B_S:=B_{\beta_S}=\Delta^{|S|+1}\times\prod_{j=1}^{r+1}C_{[\sigma(s_{j-1},\dots,s_j)]}\ \cong\ \Delta^{|S|+1}\times I^{n-1-|S|}.
\]

Every face $F$ of $B_S$ is determined by a pair $\theta=(W,\lambda)$, where $W=\{w_0<\cdots<w_p\}\subseteq\{0\}\cup(S-1)\cup\{n\}$ and $\lambda\colon\{1,\dots,n-1\}\smallsetminus(S-1)\to\{0,1,\iota\}$ is a function. Let
\[
  K_\lambda:=\{0,n\}\cup(S-1)\cup\lambda^{-1}(1),\qquad I_\lambda:=\lambda^{-1}(\iota).
\]
By the definition of $x_\theta$ in the proof of Proposition~\ref{prismatic}, $x_\theta=b_0\langle b_1|\cdots|b_p\rangle b_{p+1}$, where $b_m$ is the sub-necklace of $\sigma(w_{m-1},\dots,w_m)$ with joints $K_\lambda \cap[w_{m-1},w_m]$ and vertices $(K_\lambda \cup I_\lambda )\cap[w_{m-1},w_m]$. Thus
 $ x_\theta$ is completely determined by $(W,K_\lambda,I_\lambda)$.

The vertices of $B_S$ are determined by pairs $(g,\lambda')$ with $g \in \{0\}\cup(S-1)\cup\{n\}$ and $\lambda'\colon\{1,\dots,n-1\}\smallsetminus(S-1)\to\{0,1\}$ a function. Write $U=K_{\lambda'}$. By the formula above with $W=\{g\}$ and $I_{\lambda}=\emptyset$, the $0$-cell determined by this data corresponds to $b_0\langle\ \rangle b_1$, where $b_0$ and $b_1$ are the paths of edges of $\sigma$ through $U\cap[0,g]$ and $U\cap[g,n]$. Denote this $0$-cell by $v(U,g)$. 

We now define a map $\Psi$ from the set of pairs $(U,g)$ with $\{0,n,g\}\subseteq U\subseteq\{0,\dots,n\}$ to the defining vertices of $G_n$. Writing $A_U=\{j\in\{2,\dots,n\}:j-1\notin U\}$ we define
\[
  \Psi(U,g):=\begin{cases}c_{0,A_U},&g=0,\\ c_{1,A_U},&g=n,\\ u_{r,A_U},&g=r-1\in\{1,\dots,n-1\},\end{cases}
\]
 The map $\Psi$ determines a bijection from the vertices of $B_S$ onto those of $P_S$ compatible with the identifications $B_S\cong\Delta^{|S|+1}\times I^{n-1-|S|}\cong P_S$ (Theorem~\ref{thm:prism-comb-type}). Hence, the bijection is induced by a combinatorial isomorphism of polytopes $\psi_S \colon B_S \xrightarrow{\cong }P_S$.

We now check that the above identifications are compatible with gluings. Let $F\subseteq B_S$ and $F'\subseteq B_T$ be faces, indexed by $\theta=(W,\lambda)$ and $\theta'=(W',\lambda')$. By the proof of Proposition~\ref{prismatic}, $\chi_{\beta_S}(F)=\chi_{\beta_T}(F')$ if and only if $x_\theta=x_{\theta'}$. The set of vertices of $F$ is given by
\[
  V(F)=\{v(U,g)\,:\,g\in W,\ K_\lambda\subseteq U\subseteq K_\lambda\cup I_\lambda\}.
\]
The set $V(F)$ determines $(W,K_\lambda,I_\lambda)$, as $W$ is the set of values of $g$, $K_\lambda$ is the intersection of the sets $U$, and $K_\lambda\cup I_\lambda$ is their union. Conversely $(W,K_\lambda,I_\lambda)$ determines $V(F)$. Hence
$\chi_{\beta_S}(F)=\chi_{\beta_T}(F')$ if and only if $V(F)=V(F')$ if and only if $\Psi(V(F))=\Psi(V(F'))$.
Since $\psi_S$ and $\psi_T$ are induced by $\Psi$ on vertices, $\psi_S(F)$ and $\psi_T(F')$ have the same vertices. By Corollary~\ref{cor:prism-intersection}, the prisms $P_S$ form a polytopal complex, so faces of the prisms $P_S$ with the same vertex set coincide. Thus $\chi_{\beta_S}(F)=\chi_{\beta_T}(F')$ if and only if $\psi_S(F)=\psi_T(F')$.

It follows that the isomorphisms $\psi_S$ induce an isomorphism of the face posets of the union $\bigcup_S\chi_{\beta_S}(B_S)\subseteq\mathcal{B}_K$ and of the prismatic subdivision $G_n=\bigcup_SP_S$ from Theorem~\ref{thm:subdivision}. Therefore $\alpha(\xi)=\sum_S[e_{\beta_S}] \in C(\mathcal{B}_K)$ is the sum of the prisms of a cell complex isomorphic to this subdivision, and their union is combinatorially equivalent to $G_n$.
\end{proof}

\printbibliography

\end{document}